\documentclass[10pt]{article}
\usepackage[a4paper, left=0.9in, right=0.9in, top=0.8in, bottom=1in]{geometry}

\usepackage{url}            
\usepackage{booktabs}       
\usepackage{amsfonts}       
\usepackage{nicefrac}       
\usepackage{microtype}      
\usepackage{lipsum}
\usepackage{fancyhdr}       
\usepackage{graphicx}       
\graphicspath{{media/}}     
\usepackage[dvipsnames]{xcolor}
\RequirePackage{tgtermes}
\RequirePackage{newtxtext}
\RequirePackage{newtxmath}
\RequirePackage{bm}
\RequirePackage{endnotes}
\usepackage{setspace}
\usepackage{algpseudocode}

\usepackage{tikz}
\usepackage{algorithm,refcount,tabularx}
\usepackage{booktabs}
\usepackage{adjustbox}
\usepackage{textcomp}

\usepackage{amsthm}
\usepackage{amsmath}
\usepackage[english]{babel}
\usepackage{authblk}
\newtheorem{theorem}{Theorem}

\usepackage{natbib}
 \bibpunct[, ]{(}{)}{,}{a}{}{,}%
 \def\newblock{\ }%

\usepackage{appendix}

\makeatletter
\newcommand{\multiline}[1]{%
  \begin{tabularx}{\dimexpr\linewidth-\ALG@thistlm}[t]{@{}X@{}}
    #1
  \end{tabularx}
}
\makeatother

\newcommand\Algphase[1]{%
\vspace*{-.7\baselineskip}\Statex\hspace*{\dimexpr-\algorithmicindent-2pt\relax}\rule{\textwidth}{0.4pt}%
\vspace*{-.45\baselineskip}
\Statex\hspace*{-\algorithmicindent}\textbf{#1}%

\vspace*{-.75\baselineskip}\Statex\hspace*{\dimexpr-\algorithmicindent-2pt\relax}\rule{\textwidth}{0.4pt}%
}
\title{A \textit{K}-adaptability Approach to Proton Radiation Therapy Robust Treatment Planning }
\author[1*]{Zihang Qiu}
\author[2]{Ali Ajdari}
\author[2]{Mislav Bobi\'c}
\author[2]{Thomas Bortfeld}
\author[1]{Dick den Hertog}
\author[1]{Jannis Kurtz}
\author[3]{Hoyeon Lee}

\affil[1]{Amsterdam Business School, University of Amsterdam, Amsterdam, The Netherlands} 
\affil[2]{Department of Radiation Oncology, Massachusetts General Hospital and Harvard Medical School, Boston, MA, United States of America}
\affil[3]{Department of Diagnostic Radiology and Centre of Cancer Medicine, LKS Faculty of Medicine, University of Hong Kong, Hong Kong SAR, China}
\affil[ $*$]{z.qiu@uva.nl}

\date{}

\begin{document}

\maketitle

\begin{abstract}
Uncertainties such as setup and range errors can significantly compromise proton therapy. A discrete uncertainty set is often constructed to represent different uncertainty scenarios. A \textit{min-max} robust optimization approach is then utilized to optimize the worst-case performance of a radiation therapy plan against the uncertainty set. However, the \textit{min-max} approach can be too conservative as a single plan has to  account for the entire uncertainty set. $K$-adaptability is a novel approach to robust optimization which covers the uncertainty set with multiple ($K$) solutions, reducing the conservativeness. Solving $K$-adaptability to optimality is known to be computationally intractable. To that end, we developed a novel and efficient $K$-adaptability heuristic that iteratively clusters the scenarios based on plan-scenario performance for the proton radiation therapy planning problem. Compared to the conventional robust solution, the developed $K$-adaptability heuristic increased the worst-case CTV $D_{min}$ up to 4.52 Gy on average across five head and neck patients. The developed heuristic also demonstrated its superiority in objective value and time-efficiency compared to the competing methods we tested. 
\end{abstract}


\section{Introduction}\label{sec:Intro}
 Proton therapy damages tumor cells by irradiating the tumor regions with protons. Proton therapy treatment planning aims to determine the proton amount and irradiation trajectory to minimize damage to normal tissues while delivering enough damage to the tumor. However, proton therapy plans are particularly sensitive to uncertainties, such as setup and range errors. These uncertainties can severely compromise the plan quality, causing insufficient damage to the tumor or unnecessary damage to the normal tissues. Robust optimization is a widely used technique to manage uncertainties in proton therapy treatment planning. In the section, we provide an introduction to proton therapy treatment planning and a literature review for the  $K$-adaptability approach to robust optimization.

\subsection{Proton radiation therapy planning}
Proton therapy is a form of radiation therapy that uses proton beams instead of x-rays to treat cancers. The goal of intensity-modulated proton treatment (IMPT) planning is to optimize the weights of small ``beamlets'' so that the treatment plan delivers a sufficient dose to the tumor while minimizing damage to surrounding healthy tissues. A voxelized image, such as a CT or MRI scan, represents a patient’s anatomy \citep{image_2021,image_2022, image_2024}. From this image, a dose influence matrix is calculated based on a specific beam and beamlet configuration. This matrix describes the dose deposition from each beamlet to each voxel. By adjusting the beamlet weights, an optimal dose distribution can be achieved to meet a physician’s specifications. However, because of the large size of the dose influence matrix, proton radiation therapy planning remains computationally intensive \citep{Treatment_planning_speedup_2013,Treatment_planning_speedup_2020,Treatment_planning_speedup_2024}.\par

The min-max formulation is widely used for radiation therapy robust planning to address uncertainties, such as setup and range errors \citep{Robust_optimization_2006,Robust_optimization_2008,Robust_optimization_2009,Robust_optimization_2011}. Setup errors refer to the misalignment between the planned and actual treatment positions, which can cause the irradiation source to be offset relative to the patient. Range errors represent uncertainties in position where the proton beam stops in the patient. Both the overshoot into the normal tissues behind the tumor and undershoot where the beam stops before it reaches the end of the tumor are problematic. Together, these types of uncertainty impact the accuracy of the dose influence matrix. To reduce the effects of these uncertainties, dose influence matrices are explicitly simulated for various discrete uncertainty scenarios. Robust optimization then aims to optimize the worst-case performance regarding the discrete uncertainty scenario set. \par

\subsection{$K$-adaptability}
The $K$-adaptability approach was first introduced in \cite{First_K_adaptability} with the goal to approximate general adaptive robust optimization problems, since the latter class of problems is inherently hard to solve, especially if the decision variables are restricted to be integer. The idea of the $K$-adaptability approach is to calculate $K$ solutions in advance, before the uncertainty realizes, where $K$ is a pre-defined number. After scenario realization, the best of the $K$ calculated solutions can be chosen. The choice of the parameter $K$ affects the quality of the solutions where larger values for $K$ lead to better objective values in general. While the restriction to a small set of $K$ different reactions (instead of allowing an arbitrary number of reactions) can lead to sub-optimal solutions, having a small set of solutions can be beneficial in several applications as disaster management or in treatment planning for radiotherapy. In the latter application, the $K$-adaptability approach allows to calculate a small set of treatment plans in advance and select the best plan with respect to the daily anatomy every day. 

The $K$-adaptability approach was mainly studied for robust optimization problems. In \cite{hanasusanto2015k,subramanyam2020k} the authors derived complexity results and exact solution algorithms based on mixed-integer programming (MIP) formulations and a branch \& bound scheme which iteratively generates partitions of the uncertainty set. In \cite{ghahtarani2023double} a logic-based Benders' decomposition approach was presented to solve the $K$-adaptability problem to optimality. In \cite{kurtz2021approximation} the approximation quality depending on the parameter $K$ was analyzed and fast approximation algorithms were derived for the case where the uncertain parameters only appear in the objective function. In \cite{kurtz2024many} bounds for the number of solutions $K$ were derived which guarantee optimality for both, objective and constraint uncertainty. Machine learning methods were used in \cite{julien2024machine} to improve the performance of the branch \& bound method introduced in \cite{subramanyam2020k}. The $K$-adaptability approach was also studied in its distributionally robust (\cite{hanasusanto2016k}) and stochastic variants (\cite{buchheim2019k,malaguti2022k}), and in decision dependent information discovery setting (\cite{paradiso2022exactapproximateschemesrobust}).

When applied to single-stage problems, i.e., all decision variables have to be made at the same time, the $K$-adaptability approach reduces to the so-called \textit{min-max-min robust optimization} approach which we study in this work. The min-max-min robust problem can be formulated as 
\begin{equation}\label{eq:min-max-min}
\min_{x^1,\ldots , x^K\in \mathcal X} \max_{\xi\in \mathcal U} \min_{i=1,\ldots ,K} \ f(x^i,\xi) ,
\end{equation}
where $\mathcal X\subset \mathbb R^n$ is the feasible set, $\mathcal U$ is a given uncertainty set, and $f:\mathcal X \times \mathcal U \to \mathbb R$ a given objective function. This problem was first studied in \cite{buchheim2016min,buchheim2017min} for convex uncertainty sets. Later, improved solution algorithms and complexity results were developed for convex uncertainty sets (\cite{arslan2022min}), convex budgeted uncertainty sets (\cite{chassein2019faster}), discrete budgeted uncertainty sets (\cite{goerigk2020min}), finite uncertainty sets (\cite{buchheim2018complexity}) and for the case $K=2$ (\cite{chassein2021complexity}). The min-max-min approach was also adapted to robust regret optimization (\cite{crema2020min}) and studied for the case of smooth and strongly convex objective functions (\cite{lamperski2023min}).

It is known, that the min-max-min robust problem for small and fixed $K$ is NP-hard, even if the uncertainty set is a polyhedron (\cite{buchheim2016min}) or a finite set of scenarios (\cite{buchheim2018complexity}). Solving the problem exactly for a given $K$ is computationally extremely challenging and only possible for small-dimensional problem instances; see \cite{arslan2022min}. At the same time the amount of efficient general purpose heuristics is sparse. In \cite{buchheim2016min} a heuristic algorithm is presented for the case of convex uncertainty sets. It was later shown in \cite{kurtz2021approximation} that this algorithm achieves certain approximation guarantees which improve with increasing $K$. However, for the case of general finite uncertainty sets fast heuristics are still missing. \par

\subsection{Purpose and contribution of this paper}

The contributions of this paper are threefold. Firstly, we propose a novel heuristic for the \textit{min-max-min} type $K$-adaptability problem with a discrete uncertainty set. Secondly, we apply the proposed heuristic to the proton therapy planning problem and evaluate its efficiency. While our experiments focus on proton therapy planning, we emphasize that the heuristic is a general method applicable to problems with a discrete uncertainty set. Thirdly, we provide a proof that the $K$-adaptability problem we investigate in this study is NP-hard. \par

\par

The structure of the paper is as follows. In Section 2, we provide an overview of proton therapy treatment planning and introduce generic formulations for the nominal and robust planning optimization problems in Sections 2.1 and 2.2, respectively. In Section 2.3, we present the $K$-adaptability problem and the proof of the problem being NP-hard. Section 3 describes the novel $K$-adaptability heuristic we develop. Section 4 describes the experiment design and Section 5 reports the results of the experiments. Finally, Sections 6 and 7 provide discussion and conclusions.

\section{Proton therapy planning optimization}\label{sec:Method}
In radiation therapy, discretized scans of a patient are acquired with imaging techniques, such as computer tomography (CT) or magnetic resonance imaging (MRI), to represent their anatomy. The three dimensional scan is divided into small cubic volumes, called voxels. Each voxel is indexed and represents a unique region in the anatomy. A physician then contours the tumor target volume and the critical healthy organs (organs at risk, OAR) to be spared. The latter can be segmented automatically by an AI software, or by a dosimetrist.  Each organ or structure is represented by a set of voxels. The goal of treatment planning is to determine a set of proton beams/beamlets and their corresponding beam intensities (weights) that yield an optimal dose distribution. 
\subsection{Nominal planning optimization problem}

A common optimization problem in radiation therapy treatment planning is given by:  

\begin{equation*}\label{eq: nominal RT planning problem}
\begin{aligned}
\min_{\boldsymbol{x},\boldsymbol{d}} \quad & f({\boldsymbol{d}}) \\
\text{s.t.} \quad & g_k(\boldsymbol{d}) \leq 0, \quad k=1,\dots,m, \\ 
& \boldsymbol{D}\boldsymbol{x} = \boldsymbol{d}, \\
& \boldsymbol{x} \geq \boldsymbol{0},
\end{aligned}
\end{equation*} where $f({\boldsymbol{d}})$ denotes the planning objective function, and $g_k(\boldsymbol{d})$ denotes the $k$th planning constraint function. An example of an objective function is to maximize the minimum dose in the tumor target volume, and an example of a planning constraint is the mean dose in an organ-at-risk (OAR) volume to not exceed 30 Gy. The vector $\boldsymbol{d}$ corresponds to the dose distribution in the anatomical scenario. $\boldsymbol{D}$ is the dose-influence matrix, and $\boldsymbol{x}$ is the proton spot weights. Specifically, the dose distribution $\boldsymbol{d}$ represents the cumulative dose received by the voxels over all the beamlets, meaning that $d_i$ is the total dose delivered to voxel $i$. The dose-influence matrix $\boldsymbol{D}$ characterizes the contribution of each beamlet to the voxels. In particular, $D_{ij}$ represents the dose contribution of beamlet $j$ to voxel $i$ when the beamlet has unit intensity. Consequently, the term $D_{ij} x_j$ gives the actual dose contribution of beamlet $j$ with intensity $x_j$ to voxel $i$. \par

\subsection{Robust optimization}

Robust optimization is a standard planning approach to managing uncertainties in proton therapy. Discrete uncertainty scenario are simulated, which are the dose-influence matrices subject to the uncertainties. The robust optimization problem then optimizes for the worst-case performance. The robust optimization formulation is given by :

\begin{equation}\label{eq:robust planning problem}
\begin{aligned}
\min_{\boldsymbol{x}} \max_{\boldsymbol{D} \in \mathcal{U}} \quad &  f(\boldsymbol{D}\boldsymbol{x})\\
\textrm{s.t.} \quad &  g_{k}(\boldsymbol{D}\boldsymbol{x}) \leq 0, \quad k=1,\dots,m, \quad \boldsymbol{D} \in \mathcal{U}, \\ 
  &\boldsymbol{x}\ge\boldsymbol{0},\\
\end{aligned}
\end{equation} where $\mathcal{U}$ denotes uncertainty scenario set. We denote the dose distribution $\boldsymbol{d}$ with $\boldsymbol{D}\boldsymbol{x}$ here, because $\boldsymbol{D}$ contains  uncertainty in (\ref{eq:robust planning problem}). \par

Setup and range uncertainties are two primary uncertainty sources considered in the robust treatment planning of proton radiation therapy. The setup uncertainty represents the difference between the treatment planning and actual treatment positions of a patient, causing difference between the designated and actual placement of the proton beamlets within the patient. The range uncertainty arises from the fact there is no precise model to describe and predict the range of a proton beamlet in different tissues \citep{Range_uncertainty_2020}. Range uncertainty is crucial to a proton radiation therapy treatment because a proton releases most of its energy at the end of its path, called a Bragg peak. Misplacement of the Bragg peak compromises tumor target coverage and organs-at-risk (OARs) sparing. It is important to note that the realization of the range uncertainty of each proton beamlet can not be known exactly 
yet due to current technology limitation. However, prompt-gamma imaging \citep{Prompt_gamma_2021, Prompt_Gamma_2022, prompt_gamma_2023} and proton radiography \citep{proton_radiography_2020,proton_radiography_2020_2}  are fast-growing research topics, which aim to detect the proton range and hence reduce range uncertainty. In this study, we assume that the magnitude of the proton range can be known. To model the uncertainty scenario set  $\mathcal{U}$, dose-influence matrices are calculated explicitly subject to a combination of the setup and range uncertainties. \par

It is important to note that problem (\ref{eq:robust planning problem}) can have no feasible solution due to the first constraints $\{g_{k}(\boldsymbol{D}\boldsymbol{x}) \leq 0 :  k=1,\dots,m, \quad \boldsymbol{D} \in \mathcal{U}$\}. In proton therapy treatment planning, there exist two kinds of constraints, underdose and overdose constraints. Underdose constraints set a lower bound to the dose a structure should receive. For example, the clinical target volume (CTV) should receive at least 50 Gy of dose, $D_{min\_CTV} \ge 50$ Gy. Overdose constraints set an upper bound to the dose a structure should receive. For example, the brainstem should receive at most 57 Gy of dose, $D_{min\_Brainstem} \leq 57$ Gy. Problem (\ref{eq:robust planning problem}) always has a feasible solution if only one kind of constraints is included. When only underdose constraints exist,   $\boldsymbol{x}=\boldsymbol{a}$ is a feasible solution, where the components of $\boldsymbol{a}$ are large enough values, if $D_{ij}>0$ for all voxels $i$ belonging to the structures subject to the underdose constraints. When only overdose constraints exist, $\boldsymbol{x}=\boldsymbol{0}$ is a feasible solution. However, if both types of constraints exist, there can be no feasible solution to problem (\ref{eq:robust planning problem}). In practice, when both types of constraints are involved, a medical physicist fine-tunes the bound value of constraints to make problem (\ref{eq:robust planning problem}) feasible. In this study, we  include only overdose type constraints which ensure the feasibility of problem (\ref{eq:robust planning problem}).

\subsection{$K$-adaptable robust optimization}
As mentioned above the $K$-adaptive version of the robust optimization problem \eqref{eq:min-max-min} can be used to calculate a (small) set of $K$ different solutions which are robust against a set of scenarios. Applied to the radiation therapy treatment planning problem \eqref{eq: nominal RT planning problem} the $K$-adaptive version can be formulated as:

\begin{equation}\label{eq:min-max-min_RT}
    \begin{aligned}
        \min_{\bm{x}^1,\ldots ,\bm{x}^K\ge 0} \ \max_{\bm{D}\in \mathcal \bm{U}} \ & \min \{ f(\bm{Dx}^i): i\in\{1,\ldots ,K\}, g_k(\bm{Dx}^i)\le 0, \ \forall \ k=1,\ldots m\}. 
    \end{aligned}
\end{equation}
Problems of the latter form are known to be computationally hard to solve. However, all known complexity results in the literature are shown for binary decision vectors $\bm{x^i}$; see e.g. \cite{buchheim2017min,buchheim2018complexity,goerigk2020min}. Since problems with continuous variables are usually computationally easier to solve, the latter results do not justify the hardness of Problem \eqref{eq:min-max-min_RT}. In the following we show that the problems remains NP-hard, even for the specific structure of the radiation therapy treatment planning problem.

\begin{theorem}
    Problem \eqref{eq:min-max-min_RT} is NP-hard.
\end{theorem}
\proof{}
We reduce the decision version of the NP-hard hitting set problem to Problem \eqref{eq:min-max-min_RT}. For a given set of items $\mathcal I=\{ 1,\ldots ,n\}$, an integer $K$, and a collection of subsets of $\mathcal I$, defined as $\mathcal S=\{ S_1,\ldots ,S_T\}$, the hitting set problem answers the question if there exists a selection of at most $K$ items from $\mathcal I$, such that each set in $\mathcal S$ contains at least one of these items. The hitting set problem is known to be NP-hard; see \cite{garey2002computers}. 

We create an instance of $\eqref{eq:min-max-min_RT}$ as follows. Define $\mathcal U=\{ D^1,\ldots ,D^T \}$ where $D^l\in\mathbb R^{(n+1)\times n}$. We define $D_{jj}^l = 0$ if $j\in S^l$ and $D_{jj}^l = 1$ otherwise for all $j=1,\ldots ,n$. Furthermore we define the $n+1$-th row as the all-one vector. All other entries of $D^l$ are zero. Furthermore we define $f(d) = \sum_{i=1}^{n} d_i$ and $g_1(d) = -d_{n+1} + 1$. With the latter setup Problem \eqref{eq:min-max-min_RT} can be reformulated as
\begin{equation}\label{eq:min-max-min_RT_proof}
    \begin{aligned}
        \min_{\bm{x}^1,\ldots ,\bm{x}^k\ge 0} \ \max_{\bm{D}\in \mathcal U} \ & \min_{\substack{i=1,\ldots ,K: \\  \sum_{j=1}^{n} x_j^i\ge 1}} \ \sum_{i=1}^{n} D_{jj}x^i_j 
    \end{aligned}
\end{equation}
We show that the optimal value of the latter problem is $0$ if and only if the answer to the hitting set problem is yes. First, note that the optimal value of \eqref{eq:min-max-min_RT_proof} is non-negative, since $D\ge 0$ and $x^i\ge 0$ for all $i=1,\ldots ,K$. Assume first the answer to the hitting set problem is yes. Then there exists a set $\mathcal I_H\subseteq \mathcal I$ with $|\mathcal I_H|\le K$ such that every set $S^l$ contains at least one of the items in $\mathcal I_H$. Define the solution to Problem \eqref{eq:min-max-min_RT_proof} where $x^i_i=1$ for every $i\in\mathcal I_H$ and $x^i_j = 0$ for $j\neq i$. Then, by definition for every $D^l$, there must exist an $x^i$ s.t. $\sum_{i=1}^{n} D_{jj}x^i_j = D_{ii} x^i_i = 0$. Hence the optimal value of \eqref{eq:min-max-min_RT_proof} is zero.

For the other direction, assume the optimal value of \eqref{eq:min-max-min_RT_proof} is zero. That means that for every $D^l$ there exists an $x^i$, s.t. $\sum_{i=1}^{n} D^l_{jj}x^i_j = 0$ and $\sum_{j=1}^{n} x_j^i\ge 1$. Let $j^*_i$ be one of the indices where $x^i_{j^*_i} > 0$. Then $D^l_{j^*_ij^*_i} x^i_{j^*_i} = 0$ and by definition of $D^l$ the set $\mathcal I_H=\{ j^*_1, \ldots ,j^*_K\}$ must be a hitting set. This proves the result.

\section{Novel $K$-adaptability approach}

The $K$-adaptability problem with a discrete uncertainty set can be viewed as optimally clustering the uncertainty set $\mathcal{U}$ into $K$ scenario clusters, where each cluster is then accounted by a robust solution. The resulting $K$ robust solutions should collectively minimize the worst-case objective value against the entire uncertainty set. \par

\subsection{Overall approach}

Our $K$-adaptability heuristic takes uncertainty set $\mathcal{U}$ as input and outputs a best found set of $K$ solutions, $X_{K}^*$, in terms of worst-case performance for each $K \in \{ 1, \dots,|\mathcal{U}|\}$. The best found combination of $K$ solutions for each $K$ gives an insight into how the objective value evolves with $K$. The heuristic consists of two phases, \textit{solution generation} and \textit{solution re-distribution}. The \textit{solution generation} phase aims at generating good solutions for each $K$. This phase uses two MIP problems to cluster the scenarios and afterwards generates robust solutions for the resulting clusters. The \textit{solution re-distribution} phase aims to enhance the worst-case  objective value for each $K$ value compared to the estimates from the \textit{solution generation} phase. It performs the plan-scenario assignment again using the two MIP problems, yet with all the previously generated solutions. The objective value obtained here should be at least as high as the one obtained in \textit{solution generation} phase due to a larger solution pool.  As the crux of $K$-adaptability is scenario clustering, we first present the two MIP problems, which appear in both \textit{solution generation} and \textit{solution re-distribution} phases in Section 3.2. Then, we describe the \textit{solution generation} and \textit{solution re-distribution} phases in Sections 3.3 and 3.4. We finally provide a graphical illustration of the proposed heuristic with a four-scenarios toy example and discuss the generality of the heuristic in Section 3.5.\par

\subsection{Clustering}
Our clustering approach is to partition similar scenarios into clusters. We measure the similarity between scenarios based on the  objective value a solution can obtain on them. This is because a solution tends to yield better objective values for scenarios that closely resemble those for which the solution was created. Therefore, we can assign solutions to uncertainty scenarios in a way that optimizes the objective value with respect to the entire uncertainty set. Scenarios assigned to the same solution are considered similar to the scenarios the solution was originally created for, hence forming a cohesive cluster. \par

We use two MIP problems for the solution-scenario assignment. The first MIP problem assigns solutions to the uncertainty scenarios to optimize the worst-case objective value,

\begin{equation} \label{eq: WC_FL}
\begin{aligned}
\min_{\boldsymbol{y},w,\boldsymbol{z}} \quad &  w\\
\textrm{s.t.} \quad &  \sum_{i=1}^{|\mathcal{X}|}v_{ij}z_{ij}\leq w, \quad j = 1,\dots,|\mathcal{U}|,\\ 
&\sum_{i=1}^{|\mathcal{X}|}y_i \leq K,  \\
&\sum_{i=1}^{|\mathcal{X}|}z_{ij} = 1, \quad j = 1,\dots,|\mathcal{U}|,\\ 
  &z_{ij} \leq y_i, \quad i = 1,\dots,|\mathcal{X}|, \quad j = 1,\dots,|\mathcal{U}|,\\ 
  & y_i = \{0,1\}, \quad i = 1,\dots,|\mathcal{X}|,\\
  &z_{ij} = \{0,1\}, \quad i = 1,\dots,|\mathcal{X}|, \quad j = 1,\dots,|\mathcal{U}|,\\
\end{aligned}
\end{equation}
where $\mathcal{X}$ denotes the solution pool,  $\mathcal{U}$ denotes the uncertainty scenario set, $K$ denotes the number of solutions that can be assigned to the uncertainty scenarios, $y_i$ is a binary variable which is 1 when solution $i$ is selected for assignment, otherwise 0, $z_{ij}$ is a binary variable which is 1 when scenario $j$ is assigned to solution $i$, otherwise 0. $v_{ij}$ is the objective value of solution $i$ on scenario $j$:

\begin{equation*}
    \begin{aligned}
        v_{ij} = f(\boldsymbol{D}^j\boldsymbol{x}^i).
    \end{aligned}
\end{equation*} $\boldsymbol{v}$ is obtained by calculating the objective value of each solution $i$ from $\mathcal{X}$ and each uncertainty scenario $j$ from $\mathcal{U}$.
Multiple assignments with the same optimal worst-case objective value can exist yet with different scenario-average objective values. Therefore, a second MIP problem is needed to calculate, given the optimal worst-case performance, the assignment with the best average objective value:

\begin{equation} \label{eq: avg_FL}
\begin{aligned}
\min_{\boldsymbol{y},\boldsymbol{z}} \quad &\sum_{i=1}^{|\mathcal{X}|}\sum_{j=1}^{|\mathcal{U}|}v_{ij}z_{ij}\\
\textrm{s.t.} \quad 
&\sum_{i=1}^{|\mathcal{X}|}v_{ij}z_{ij}\leq w^{*}, \quad j = 1,\dots,|\mathcal{U}|,\\ 
&\sum_{i=1}^{|\mathcal{X}|}y_i \leq K,  \\
&\sum_{i=1}^{|\mathcal{X}|}z_{ij} = 1, \quad j = 1,\dots,|\mathcal{U}|,\\ 
  &z_{ij} \leq y_i, \quad i = 1,\dots,|\mathcal{X}|, \quad j = 1,\dots,|\mathcal{U}|,\\ 
  & y_i = \{0,1\}, \quad i = 1,\dots,|\mathcal{X}|,\\
  &z_{ij} = \{0,1\}, \quad i = 1,\dots,|\mathcal{X}|, \quad j = 1,\dots,|\mathcal{U}|,\\
\end{aligned}
\end{equation} where $w^*$ is the optimal worst-case objective value from problem (\ref{eq: WC_FL}). The assignment from the problem (\ref{eq: avg_FL}) is optimal in both worst-case and scenario-average objective values with respect to the solution pool for assignment and uncertainty scenarios, achieving \textit{Pareto robust optimality} \citep{Pareto_robust_optimality_2014, Pareto_robust_optimality_2022}.

\subsection{Solution generation}

The \textit{solution generation} phase follows an iterative approach to progressively enhance the objective value. First, a global solution pool is constructed by generating an optimal solution for each uncertainty scenario. Then, a loop to generate more solutions executes in a descending order in $K$. For each $K$ value, at each iteration, problems (\ref{eq: WC_FL}) and (\ref{eq: avg_FL}) select $K$ solutions from the global solution pool. The uncertainty scenarios assigned to these selected solutions form $K$ clusters. A robust solution is generated for each resulting cluster and added to the pool. As the global solution pool expands, Problems (\ref{eq: WC_FL}) and (\ref{eq: avg_FL}) might find better combinations of $K$ solutions in subsequent iterations, leading to an improved worst-case objective value. The iterative process stops when  $K$ scenario clusters $\boldsymbol{C}_{K}$ that have already appeared previously appear again. The algorithm terminates in finite time because there are only finitely many ways to partition the $N$ scenarios into $M$ non-empty sets.

\subsection{Solution re-distribution}
The \textit{solution re-distribution} phase is devised to improve the worst-case objective value. A best combination of $K$ solutions is selected from \textit{all} the solutions generated in the \textit{solution generation} phase by solving Problems (\ref{eq: WC_FL}) and (\ref{eq: avg_FL}). Since the solution pool in this phase includes solutions generated for all $K$ values, the objective value will be at least as high as that estimated in the \textit{solution generation} phase, where only solutions generated up to the current $K$ value are available for assignment.

\subsection{Graphical illustration and generality of the $K$-adaptability heuristic}

\begin{figure}[!hbt]
    
    \centering
    \label{fig:K-adaptability_schematic}
    \caption{Graphical illustration of the proposed $K$-adaptability heuristic with four uncertainty scenarios. A blue square denotes an uncertainty scenario, an orange triangle denotes a treatment plan, $w_{ij}*$ and $\boldsymbol{C}_{ij}^*$ denote the worst-case objective value and scenario clusters obtained for $K=i$ at iteration $j$. MIP and RO are acronyms of mixed-integer programming problem and robust optimization.}
    \includegraphics[width=1.0\textwidth]{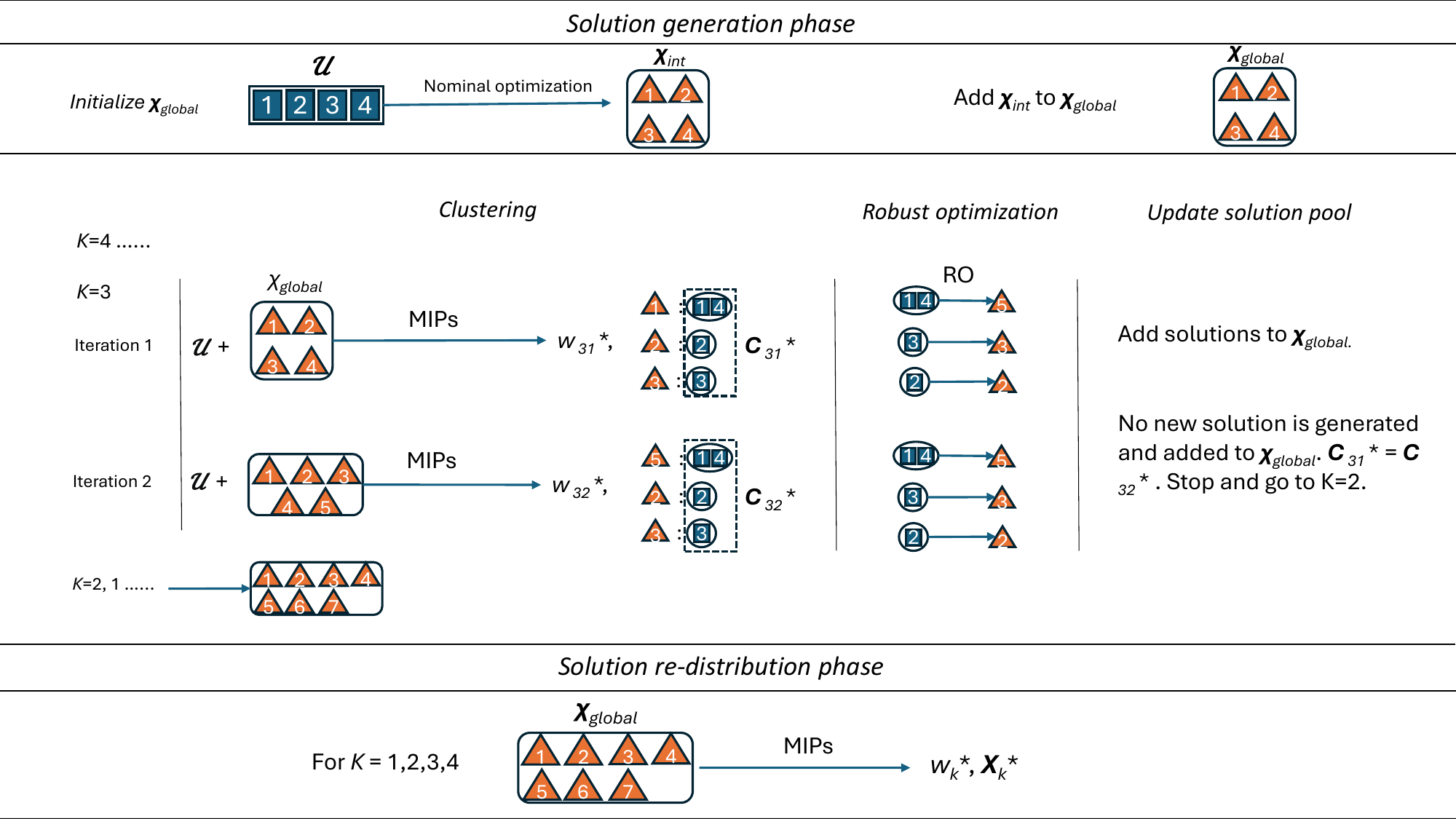}

\end{figure}

Algorithm \ref{alg:K_adaptability_heuristic} describes our $K$-adaptability heuristic. Figure \ref{fig:K-adaptability_schematic} exemplifies the heuristic with a four-scenarios toy problem. The \textit{solution generation} phase begins with generating a solution for each scenario in $\mathcal{U}= \{\boldsymbol{D}_1, \boldsymbol{D}_2,\boldsymbol{D}_3,\boldsymbol{D}_4\}$ forming an initial solution pool $\{\boldsymbol{x}_1,\boldsymbol{x}_2,\boldsymbol{x}_3,\boldsymbol{x}_4\}$. The initial solution pool is then added to the global solution pool. We skip visualizing the iterative solution generation process at $K=4$ because each initial solution will be assigned to their corresponding scenario at $K=4$. Consequently, no new solution will be generated at $K=4$. At $K=3$, the two MIP problems select 3 solutions from the global solution pool $\mathcal{X}_{global}$ that maximize the worst-case and scenario-average objective value against uncertainty set $\mathcal{U}$. In the first iteration, $\boldsymbol{x}_1$ is assigned to the $\boldsymbol{D}_1$ and $\boldsymbol{D}_4$, $\boldsymbol{x}_2$ is assigned to $\boldsymbol{D}_2$, and  $\boldsymbol{x}_3$ is assigned to $\{\boldsymbol{D}_3\}$. This results in worst-case scenario objective value $w_{31}^*$ and scenario clusters $\boldsymbol{C}_{31}*= \{\{\boldsymbol{D}_1,\boldsymbol{D}_4\}, \{\boldsymbol{D}_2\}, \{\boldsymbol{D}_3\}\}$. Then, an optimal solution is generated for each cluster in $\boldsymbol{C}_{31}^*$ and added to $\mathcal{X}_{global}$. In the second iteration, the two MIPs assign $\boldsymbol{x}_5$, which is optimal with respect to $\{\boldsymbol{D}_1,\boldsymbol{D}_4\}$, to $\boldsymbol{D}_1$ and $\boldsymbol{D}_4$, $\boldsymbol{x}_2$  to $\boldsymbol{D}_2$, and  $\boldsymbol{x}_3$  to $\boldsymbol{D}_3$. No new solution is added to $\mathcal{X}_{global}$ because $\boldsymbol{C}_{32}^* = \boldsymbol{C}_{31}^*$, and the heuristic hence stops for $K=3$ and proceeds with generating solutions for $K=2$.  After generating solutions for $K=2$ and $K=1$, $\mathcal{X}_{global} $ expands to $\{\boldsymbol{x}_1,\boldsymbol{x}_2,\boldsymbol{x}_3,\boldsymbol{x}_4,\boldsymbol{x}_5,\boldsymbol{x}_6, \boldsymbol{x}_7\}$ at the end of the \textit{solution generation} phase. In the \textit{solution re-distribution} phase, the two MIP problems select the optimal combination of $K$ solutions $X_{K}^*$ for $K=1,\dots,4$ from the complete global solution pool $\mathcal{X}_{global}=\{\boldsymbol{x}_1,\boldsymbol{x}_2,\boldsymbol{x}_3,\boldsymbol{x}_4,\boldsymbol{x}_5,\boldsymbol{x}_6, \boldsymbol{x}_7\}$.

While  the $K$-adaptability heuristic was originally developed for the robust treatment planning problem in proton therapy, it is essentially a generic method for robust optimization with discrete uncertainty sets. Therefore, the developed heuristic is also applicable to the same type of problem in other domains, such as logistics \citep{Logistic_discrete_2024, Logistic_discrete_2024_2} and knapsack \citep{knapsack_discrete_1998, knapsack_discrete_2006} problems.


\begin{algorithm}[!hbt]
\caption{$K$-adaptability heuristic}\label{alg:K_adaptability_heuristic}
\vskip6pt
\textbf{Input:} Uncertainty scenario set $\mathcal{U}$ = $\{\boldsymbol{D}^1, \dots, \boldsymbol{D}^{T}\}$.\\
\textbf{Output:} Best found combination of $K$ solutions with respect to $\mathcal{U}$ in terms of worst-case objective value $w^*_K$: $X_{K}^*$ = $\{\boldsymbol{x}_1^{(K)}, \dots, \boldsymbol{x}_K^{(K)}\}$, for $K$ = 1,\dots,$|\mathcal{U}|$.
\begin{algorithmic}[1]
\Algphase{Phase 1 - Solution generation}
\State Initialize $K \gets |\mathcal{U}|$ and global solution pool $\mathcal{X}_{global} \gets \emptyset$.
\State Generate an optimal solution for each scenario in $\mathcal{U}$ to form the initialization solution set $\mathcal{X}_{int}$.
\State $\mathcal{X}_{global} \gets  \mathcal{X}_{global} \cup \mathcal{X}_{int}$ .
\While{$ K\ge1$}
    \State Initialize scenario clusters pool $\boldsymbol{\mathcal{P}} \gets \emptyset$ and $\boldsymbol{C}_{K} \gets \emptyset$
    \Repeat
        \State $\boldsymbol{\mathcal{P}}\gets \boldsymbol{\mathcal{P}} \cup \boldsymbol{C}_{K}$
        \State Solve problem (\ref{eq: WC_FL}) on  $\mathcal{X}_{global}$. Obtain   $w^{*}$.
        \State Solve problem (\ref{eq: avg_FL}) on  $\mathcal{X}_{global}$ subject to $w^{*}$. Obtain $\boldsymbol{C}_{K}$. 
         \State Solve problem (\ref{eq:robust planning problem}) for each scenario cluster in $\boldsymbol{C}_{K}$. Obtain solutions   $\mathcal{X}_{R}$.
        \State  $\mathcal{X}_{global} \gets  \mathcal{X}_{global} \cup \mathcal{X}_{R}$.
    \Until $\boldsymbol{C}_{K} \in \boldsymbol{\mathcal{P}}$
   \State $K \gets  K - 1$ 
\EndWhile 
\Algphase{Phase 2 - Solution re-distribution}
\State Set $K \gets |\mathcal{U}|$
\While{$ K\ge1$}
    \State Solve problem (\ref{eq: WC_FL}) on  $\mathcal{X}_{global}$. Obtain $w^*_K$.
    \State Solve problem (\ref{eq: avg_FL}) on  $\mathcal{X}_{global}$ subject to $w^{*}_K$. Obtain $X_{K}^*$.
    \State $K \gets  K - 1$
\EndWhile
\end{algorithmic}
\end{algorithm}

\section{Experiment setup}\label{sec:Experiments}
This section describes the information on the data used for the computational experiment, methods we compare the $K$-adaptability heuristic with, and performance evaluation metrics. 
\subsection{Data and software}

The patient cohort includes five head and neck (H\&N) patients obtained from Massachusetts General Hospital. The CT images of the patients were acquired by a wide-bore GE scanner
(General Electric Medical Systems, Milwaukee, WI). The structure contours were manually delineated by experienced medical physicists at Massachusetts General Hospital. Figure \ref{tab:Patient_info} contains the information on the voxel number of the region of interest (ROI). ROI consists of structures considered in treatment planning. 

\begin{table}[hbt!]
\caption{\label{tab:Patient_info} Patient information. ROI denotes the region of interest.}
\centering
\begin{tabular}{lrr}\toprule

Case   &ROI voxel number  & Proton beamlet number \\ \midrule

 1    & 47,027 &  7,414      \\ 
 
 2       &   57,405   & 9,201                 \\ 
  3      &    29,457  & 7,747              \\
 4         &   47,797    & 6,245               \\ 
  
 5         &   63,731   & 8,120              
         \\ \bottomrule\end{tabular}

\end{table}

A dose of 57 Gy is prescribed to the CTV. The treatment planning objective function is to maximize the minimum dose within the CTV. The treatment planning constraints are based on the QUANTEC report recommendation for H\&N cancer \citep{quantec_2010, Constraint_functions}. Table \ref{tab:planning_constraints} exhibits the specification of the used planning constraint functions. The CTV 5 mm rind denotes the volume formed by isotropically expanding the CTV by 5 mm. A three treatment beam configuration was used in the study. The gantry angles of the beams are 180, 300, and 300\textdegree, respectively. The couch angles of the beams are 0, 180, and 0\textdegree, respectively.

\begin{table}[hbt!]
\caption{\label{tab:planning_constraints} The table of planning constraint functions.}
\centering
\begin{tabular}{llr}\toprule

Structure & Constraint function      & Bound value (Gy)  \\ \midrule

 CTV &    Maximum dose    &  59.85      \\ 
 
 CTV 5 mm rind &    Maximum dose               & 57.00                 \\ 
  Brainstem &      Maximum dose            & 54.00               \\
 Constrictors  &   Mean dose                & 25.00               \\ 
  
 Larynx &  Mean dose                & 40.00              \\ 
 Left parotid&   Mean dose                        & 26.00                   \\ 
 Right parotid&   Mean dose                        & 26.00                   \\ 
 Spinal cord&   Maximum dose                        & 45.00                   \\ 
         \\ \bottomrule\end{tabular}

\end{table}
Monte Carlo dose calculation engine MOQUI \citep{MOQUI} was used for dose influence matrix calculation and setup and range uncertainty application. For each patient, 57 uncertainty scenarios were generated. The uncertainty scenarios are combinations of a setup uncertainty and a range uncertainty. The magnitude of the setup uncertainty is 3 mm. A total of 19 directions of the setup uncertainty shift are considered in the study: one with no shift, six shifts along one of the x, -x, y, -y, z, -z directions, and twelve shifts along the combination of two of the x, -x, y, -y, z, -z directions. Three range errors, +3, -3, and 0\%, were employed in the study because 3 \% is a common estimate of range uncertainty \citep{Range_uncertainty_2020, Range_uncertainty_3_per_2023}. This study always includes the nominal scenario, i.e. no shift and range uncertainties, in the proton therapy treatment planning optimization problem (\ref{eq:robust planning problem}). This is because a robust treatment plan is always expected to control the dose distribution on the nominal scenario by the clinical practice.  \par

Gurobi 11.0.3 (Gurobi Optimization, Inc., Houston, TX) was used to solve the MIP problems $(\ref{eq: WC_FL})$, $(\ref{eq: avg_FL})$, and $(\ref{eq: $K$-medoid})$. Nymph 2023.11.09 \citep{Nymph} was used to solve the radiation therapy treatment plan robust optimization problem, i.e. problem $(\ref{eq: $K$-medoid})$. Nymph is the clinical proton therapy treatment plan optimizer at Massachusetts General Hospital.  An AMD Genoa CPU with 48 cores and 49.75 GB RAM were used to solve the optimization problems above.\par

\subsection{Competing methods}
We compare our $K$-adaptability heuristic to three competing clustering methods to demonstrate its effectiveness. The three competing methods are the local solution pool (LSP) and ascending order solution generation (AOSG) variants of the proposed $K$-adaptability heuristic, and the $K$-medoids method. The LSP and AOSG variants are simple variants of our heuristic developed by us.\par
The LSP variant employs a local solution pool in the \textit{solution generation} phase in contrast to the global solution pool employed by the $K$-adaptability heuristic. The global solution pool accumulates the previously generated solutions.  On the other hand, the local solution pool of the LSP variant always starts with the ${|\mathcal{U}|}$ initial optimal solutions generated for the uncertainty scenarios. Solutions generated for the previous $K$ values do not enter the local solution pool. Algorithm \ref{alg:LSP_heuristic_variant} describes the LSP variant.\par
The AOSG variant differs in the order of solution generation. The $K$-adaptability heuristic generates solutions in a descending order in $K$, namely from $|\mathcal{U}|$ to 1. In contrast, the AOSG variant generates solutions in an ascending order in $K$, namely from 1 to $|\mathcal{U}|$. Algorithm \ref{alg:AOSG_heuristic_variant} describes the AOSG variant.\par

Unlike the $K$-adaptability heuristic and its variants that cluster based on the performance of a given set of solutions with respect to the uncertainty scenarios, the $K$-medoids method clusters directly based on the uncertainty scenarios. $K$-medoids minimizes the sum of the distance between the uncertainty scenarios and one of the $K$ centroids \citep{K_medoids_2008,K_medoids_2010}. The $K$ centroids are selected from the uncertainty scenarios. For the algorithm used to solve the $K$-medoids method, see appendix \ref{Appendix:k-medoids}.

\section{Experiment results}

In this section, we present the numerical results of applying the $K$-adaptability approach to a proton therapy robust treatment planning problem. The comparison with conventional robust optimization illustrates the improvement that can be obtained for the current clinical practice by employing the $K$-adaptability approach. The comparison with the competing clustering methods demonstrates the effectiveness of the developed $K$-adaptability heuristic. For dose distribution, dose-volume histogram (DVH) and objective plots of individual cases, see Appendix \ref{Appendix:individual result}.

\subsection{$K$-adaptability vs conventional robust optimization}

We compare the $K$-adaptability heuristic with conventional robust optimization to assess the dosimetric improvements over current clinical standard practice. Specifically, we examine the worst-case objective values for $K$=1, 2, 3, 4, 5 and 57. At $K$=1, the worst-case objective value corresponds to that of conventional robust optimization. At $K$= 2, 3, 4, and 5, the worst-case objective value corresponds to the performance that can be obtained by employing the $K$-adaptability heuristic with a reasonable number of plans. Creating a treatment plan can cost significant amount of time and human resource in quality assurance (QA). We therefore deem preparing up to five plans reasonable. At $K$=57, the worst-case objective value represents the best possible performance achievable with the $K$-adaptability approach, because each scenario gets its nominal plan assigned in this case. The difference between these two values quantifies the improvement in objective value that can be achieved by replacing conventional robust optimization with the $K$-adaptability approach.\par

\begin{figure}[hbt!]
    
    \centering
    \caption{Worst-case dose distribution of case four. Upper left and right figures are the dose distributions of the conventional robust plan and $K$-adaptability ($K$=57), respectively. The grey area denotes the CTV. The bottom figure is the CTV DVH of the two plans with respect to the worst-case scenario.}
    \includegraphics[width=1.0\textwidth]{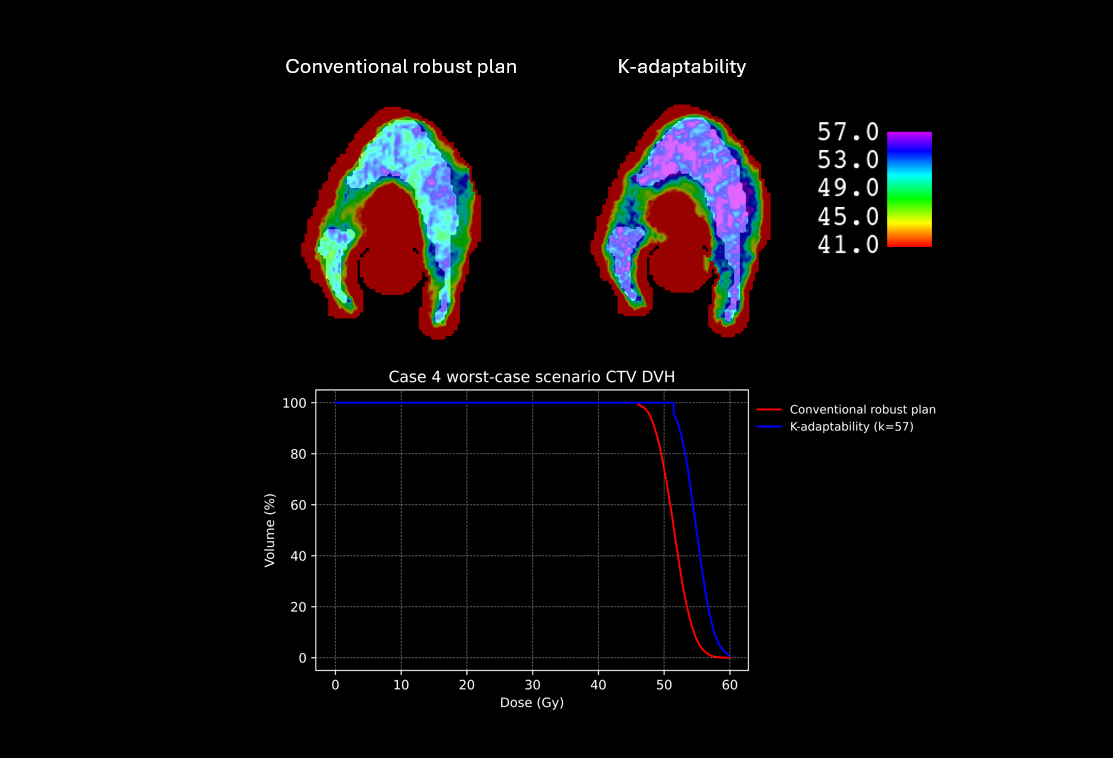}
    \label{fig:case_4_dvh}

\end{figure}

Table \ref{tab:WC_CTV_improvement} presents the improvement by employing the $K$-adaptability heuristic from the conventional robust solution ($K$=1). On average, the minimum dose to the clinical target volume (CTV), $D_{min}$, increased by 4.52 Gy across the five cases, with a maximum improvement of 5.38 Gy and a minimum of 3.75 Gy. Figure \ref{fig:case_4_dvh} illustrates the worst-case dose distribution obtained with the conventional robust optimization and the $K$-adaptability approach ($K$=57) in case 4. In the conventional robust plan, underdosage is observed in the CTV, which receives an average of 51.54 Gy. In contrast, the $K$-adaptability plan significantly enhances target coverage, delivering an average of 55.00 Gy to the CTV.

\begin{table}[hbt!]
\caption{\label{tab:WC_CTV_improvement} Worst-case CTV $D_{min}$ improvement (Gy) from the conventional robust solution ($K$=1).}
\centering
\begin{tabular}{lrrrrr}\toprule

Case  & $K$ = 2&$K$=3 &$K$=4 &$K$=5 & $K$=57 \\ \midrule

 1   &0.93 &1.39 &2.25 &2.38 & 4.16       \\ 
 
 2    & 0.99&1.93 &2.42 &2.42   &   4.62                    \\ 
 3   &1.07 &1.89 &2.12 &2.78   &    4.67             \\
 4    & 1.35&1.94 &2.75 &2.75     &   5.38                 \\ 
  
 5     & 0.86&1.39 &2.18 &2.50    &   3.75                \\ 
 Overall&$1.04\pm0.17$ & $1.71\pm0.26$& $2.34\pm0.23$ &$ 2.57\pm0.17$& $4.52\pm0.54$
         \\ \bottomrule\end{tabular}

\end{table}



\subsection{Comparison with competing clustering methods}

We evaluate the performance of the $K$-adaptability heuristics in two aspects: objective value and computational time. Since solving the $K$-adaptability problem exactly for proton therapy treatment planning is computationally intractable, we lack precise worst-case objective values for each $K$. Therefore, we compare heuristics relative to each other, with the heuristic yielding the highest objective value considered the best-performing one. For objective value, we analyze: 1. The sum of the worst-case objective value from $K$=1 to 10, which reflects heuristic performance at lower $K$ values when the objective value fluctuates more. 2. The saturation $K$ value, defined as the smallest $K$ for which the objective value equals the best possible value at $K$=57. A lower saturation $K$ indicates better performance. The saturation $K$ value can also be interpreted as the minimum number of treatment plans required to achieve the maximum worst-case performance under the uncertainty set. For computational time, we report the total runtime of each heuristic, including data loading and saving, matrix calculations, and optimization.\par

\begin{figure}[hbt!]
    
    \centering
    \caption{Worst-case dose distribution of case four. The plot exhibits the worst-case $D_{min}$ evaluated by the heuristics at different $K$ values. The $K$-adaptability heuristic shows superior assessment of the saturation $K$ value and the objective value sum from $K$ = 1 to 10.}
    \includegraphics[width=1.0\textwidth]{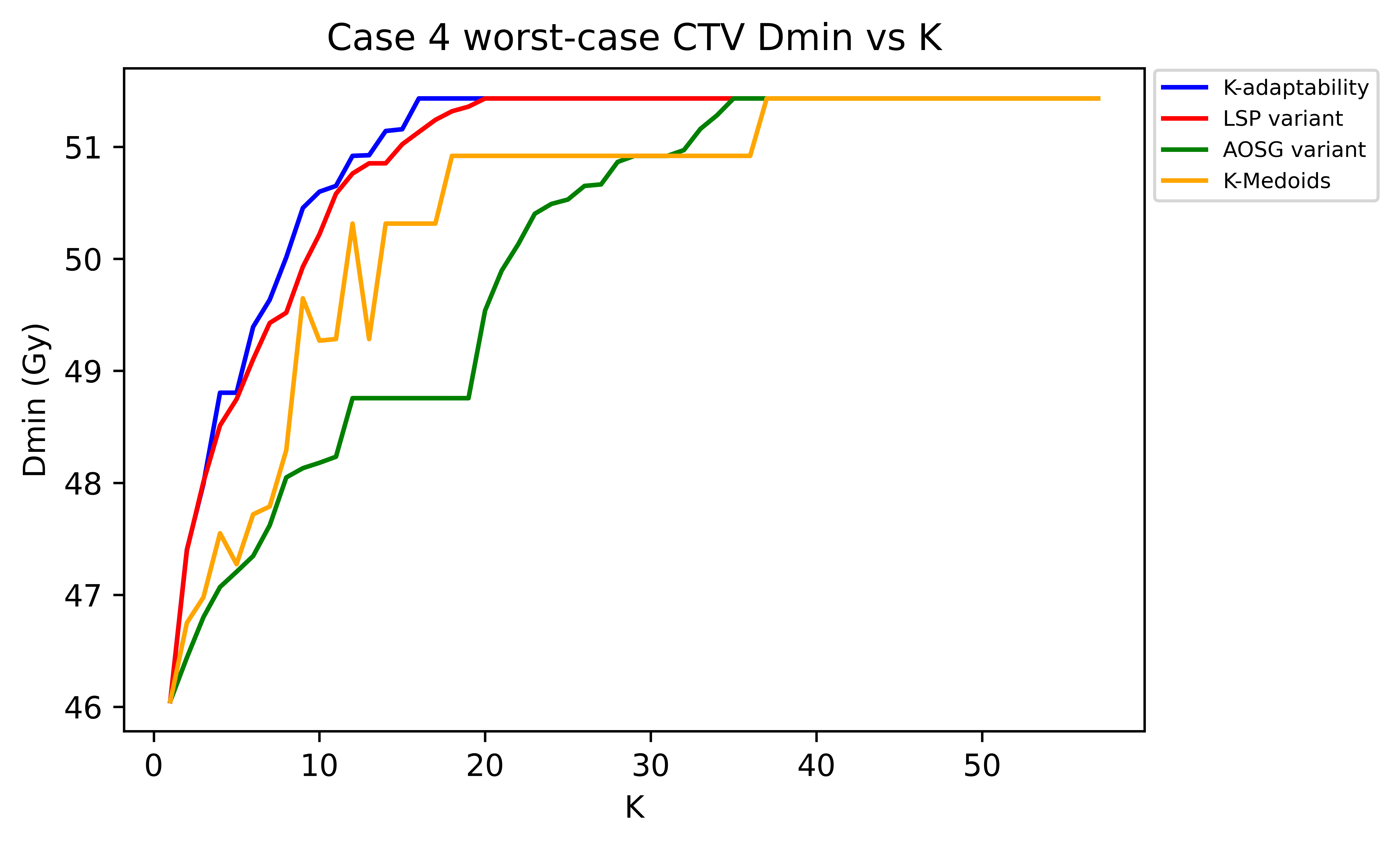}
    \label{fig:case_4_worst_case_objective_vs_K}

\end{figure}

Table \ref{tab:Heuristic_objective_value} displays the objective value metrics evaluated by the methods. The developed $K$-adaptability heuristic and the LSP variant had similar objective value sum until $K$ = 10 on average over the five cases, being 48.34 and 48.26, respectively. The $K$-adaptability heuristic had a slight edge in saturation $K$ value compared to the LSP variant, with the average values being 15.00 and 16.60, respectively. The AOSG variant and $K$-medoids methods are inferior in terms of objective value. Figure \ref{fig:case_4_worst_case_objective_vs_K} exhibits the worst-case objective value vs $K$ plot of case four. Specifically, $K$-medoids has average $K$ saturation value being 49.60. 

\begin{table}[hbt!]
\centering
\caption{\label{tab:Heuristic_objective_value} Objective value metrics evaluated by the heuristics over five cases. Sum is the objective value sum over $K$ = 1 to 10 in Gy and saturation is the saturation $K$ value.}
\begin{adjustbox}{width=1.0\textwidth}

\begin{tabular}{lrcrcrcrcrcrcrr}
\hline
& \multicolumn{2}{c}{Case 1} & \multicolumn{2}{c}{Case 2}& \multicolumn{2}{c}{Case 3}& \multicolumn{2}{c}{Case 4}& \multicolumn{2}{c}{Case 5}& \multicolumn{2}{c}{Average}\\ \cmidrule(r){2-3}  \cmidrule(r){4-5} \cmidrule(r){6-7} \cmidrule(r){8-9} \cmidrule(r){10-11} \cmidrule(r){12-13}

        Heuristic & Sum & Saturation   & Sum& Saturation  & Sum & Saturation  & Sum& Saturation  & Sum& Saturation & Sum& Saturation \\ 
        \midrule
        $K$-adaptability &45.78 & 10  & 50.55 &21 & 49.23& 18& 48.91&16 & 47.23&10 & 48.34&15.0\\
        LSP variant  &45.74 & 12  &  50.45 & 22&49.16 & 19& 48.69&20& 47.25& 10& 48.26&16.6\\
        AOSG variant&44.42 & 20  & 49.11 & 38&47.75 & 38&47.29 &35&45.61 &18 & 46.84&29.8\\
        $K$-medoids&45.12 & 56 & 49.97  & 57&48.62 & 51& 47.73&37& 46.62& 47&47.61 &49.6\\
        \hline

\end{tabular}
\label{tab:heuristics}
\end{adjustbox}
\end{table}



\par

Table \ref{tab:Heuristic_timing} presents the number of Nymph optimization runs and the runtime for each heuristic. The AOSG variant and the $K$-medoids method are the most time-efficient, with average runtimes of 45,256 and 46,592 seconds, respectively. The $K$-adaptability heuristic follows, with an average runtime of 65,314 seconds. The LSP variant is the slowest, averaging 83,766 seconds. Notably, the $K$-adaptability heuristic is approximately 28\% faster than the LSP variant with 10\% fewer Nymph runs.

\begin{table}[hbt!]
\centering
\caption{\label{tab:Heuristic_timing} Table of the number of nymph optimization run and runtime (s) of the heuristics in the five cases.}
\begin{adjustbox}{width=1.0\textwidth}

\begin{tabular}{lrrrrrrrrrrrrrr}
\hline
& \multicolumn{2}{c}{Case 1} & \multicolumn{2}{c}{Case 2}& \multicolumn{2}{c}{Case 3}& \multicolumn{2}{c}{Case 4}& \multicolumn{2}{c}{Case 5}& \multicolumn{2}{c}{Average}\\ \cmidrule(r){2-3} \cmidrule(r){4-5} \cmidrule(r){6-7} \cmidrule(r){8-9} \cmidrule(r){10-11} \cmidrule(r){12-13}

        Heuristic & \# Nymph opt. & Runtime & \# Nymph opt. & Runtime & \# Nymph opt. & Runtime & \# Nymph opt. & Runtime & \# Nymph opt. & Runtime & \# Nymph opt. & Runtime \\ 
\midrule
$K$-adaptability  & 223 & 72,396  & 229 & 64,343  & 203 & 55,326  & 218 & 48,975  & 195 & 85,530 & 213.6& 65,314  \\

LSP variant       & 232 & 89,364       & 262 & 86,149      & 233 & 72,335       & 231 & 62,187      & 221 & 108,794       & 235.8 &    83,766    \\
AOSG variant       & 136 & 43,446  & 114 &  39,795  & 123 & 43,028 & 140 & 38,230 & 135 & 61,782  & 129.6 & 45,256  \\

$K$-medoids        &168 & 47,904 & 148 &40,883 & 162 & 40,268 & 171 & 36,660 & 169 &  67,244 & 163.6 & 46,592 \\
\hline
\end{tabular}
\end{adjustbox}
\end{table}

\section{Discussion}

This study compares the worst-case CTV $D_{\min}$ achieved using conventional robust optimization and the $K$-adaptability approach. The results demonstrate a substantial improvement of 4.52 Gy on average with the $K$-adaptability approach, highlighting its dosimetric benefits. These findings support the potential advantages of multi-plan treatment, where multiple treatment plans are available for daily selection. To establish multi-plan treatment as a viable clinical approach, technical and logistical challenges, such as the time and labor costs associated with quality assurance (QA), must be addressed. The development of automated QA tools could play a crucial role in facilitating its implementation.  \par

Our $K$-adaptability approach estimated that, on average, 15 treatment plans are needed to optimize the worst-case CTV $D_{\min}$ across a 57-scenario uncertainty set, based on the five cases analyzed in the study. The saturation plan number is approximately one-fourth of the uncertainty set size, indicating that the maximum worst-case performance can be achieved with a plan set significantly smaller than the full uncertainty set. \par  

The developed $K$-adaptability heuristic showed superior objective value compared to the common $K$-medoids method. There are two key reasons why the $K$-adaptability heuristic outperforms the $K$-medoids method in this aspect.  First, the $K$-medoids method clusters scenarios based on the $L$-2 norm between dose-influence matrices, whereas the $K$-adaptability heuristic clusters scenarios based on plan-on-scenario performance—specifically, the objective value—which is a more directly relevant measure. Second, the $K$-medoids method minimizes the sum of $L$-2 norms between the dose-influence matrices and their respective centroids, as shown in Equation~(\ref{eq: $K$-medoid}). In other words, it optimizes for average scenario performance, whereas the $K$-adaptability approach explicitly accounts for worst-case performance, as defined by Equation (\ref{eq: WC_FL}).  
 \par

\begin{figure}[hbt!]
    
    \centering
    \caption{Size of the cluster that contains the worst-case scenario vs K.}
    \includegraphics[width=0.8\textwidth]{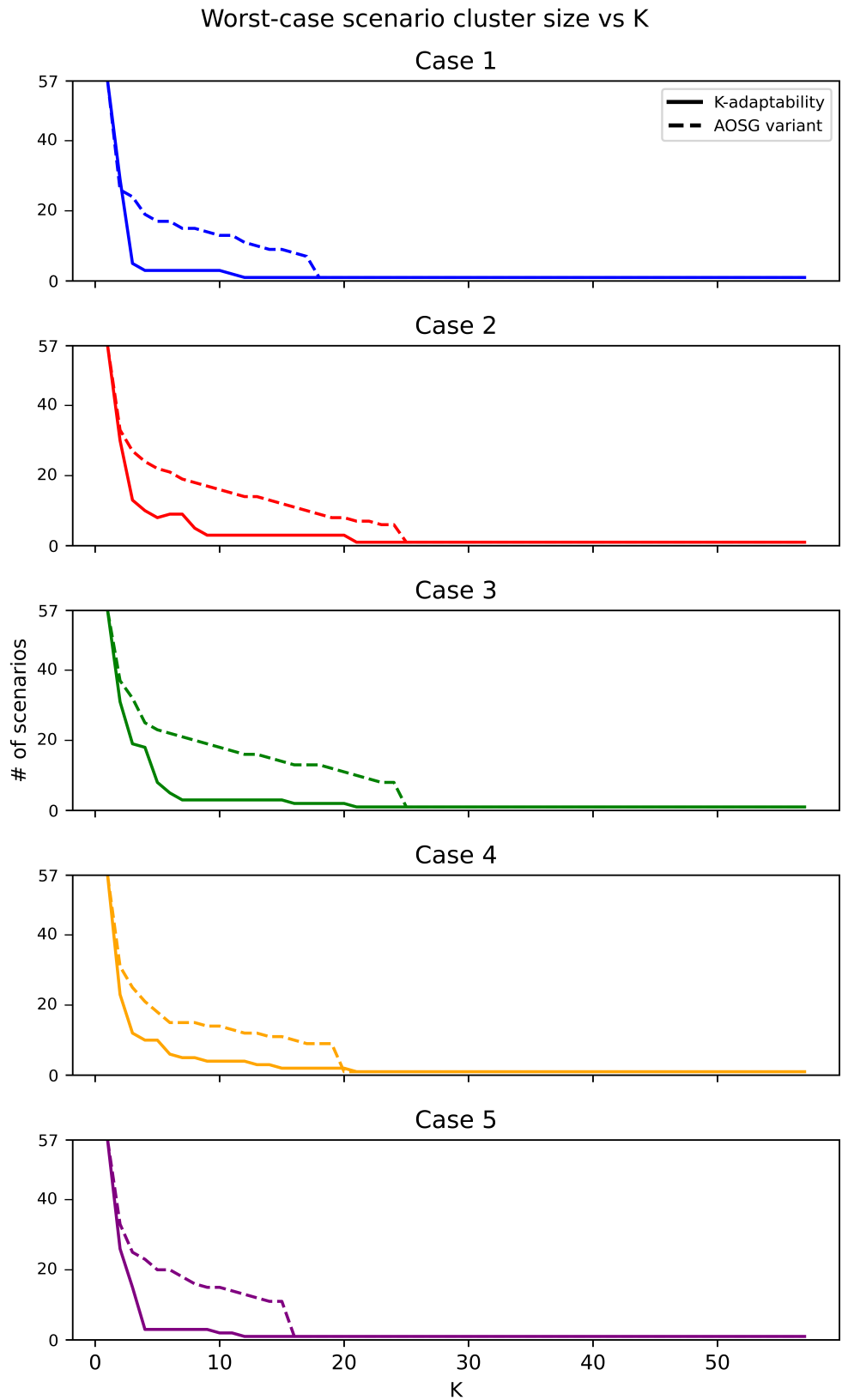}
    \label{fig:Worst-case-cluster-size}

\end{figure}

The $K$-adaptability heuristic achieved better objective values than the AOSG variant, particularly at lower $K$ values. This suggests that the order in which solutions are generated affects the performance of the heuristic. The difference in objective values is attributed to the different sizes of the cluster containing the worst-case caused by the \textit{solution generation} order in $K$. A larger scenario cluster is likely to lead to a more conservative solution that is worse in the objective value.  Figure \ref{fig:Worst-case-cluster-size} shows the size of the cluster containing the worst-case scenario of the $K$-adaptability heuristic and the AOSG variant in the \textit{solution generation} phase. The worst-case scenario cluster is larger for the AOSG variant than for the $K$-adaptability heuristic at lower $K$ values. The formation of a larger worst-case scenario cluster at lower $K$ for the AOSG variant is likely due to the dominance of robust solutions. When solutions are generated in ascending order, at $K = 2$, the solution pool consists only of the initial solutions and the robust solution from $K = 1$. As the robust solution is optimized for all scenarios, it gets assigned to most scenarios, including the worst-case scenario, forming a large cluster. In contrast, when solutions are generated in descending order, the pool at $K = 2$ includes solutions from higher $K$, namely from $|\mathcal{U}|$ to 3. The larger solution pool offers more assignment options, leading to better clustering. Additionally, these solutions are less conservative than the robust solution at $K = 1$, allowing for a more even distribution of scenarios. In the end, the larger worst-case cluster of the AOSG variant leads to a more conservative solution, resulting in a worse objective value compared to the $K$-adaptability heuristic.



The comparison with the LSP variant suggests that a global solution pool enhances heuristic time efficiency while yielding slightly better objective estimates. The $K$-adaptability heuristic is 28\% faster than the LSP variant and, on average, requires 10\% fewer Nymph runs. The higher number of Nymph runs in the LSP variant may come from the fact that its local solution pool only starts with the $|\mathcal{U}|$ initial solutions for selection. Since these initial solutions do not dominate each other, the uncertainty scenarios are evenly distributed into $K$ clusters, leading to the generation of more unique clusters. Consequently, Nymph is called to generate solutions for these unique clusters. In contrast, the global solution pool of the $K$-adaptability heuristic retains previously generated solutions. These solutions are likely to be reassigned to the same scenarios at the next $K$ value. This results in repeated clusters with existing solutions, for which Nymph is not called for solution generation. 

A future research direction would be comparing the $K$-adaptability approach to the conventional robust optimization approach on patients acquired temporal data, i.e, daily images. In the study, we evaluate $K$-adaptation and conventional robust optimization approaches on the training uncertainty set where the uncertainty scenarios are created with setup and range uncertainties. In radiation therapy, the setup and range uncertainties are used as a surrogate for anatomical changes because anatomical changes are hard to predict. The resulting robust solution in fact aims at accounting for the anatomical changes to occur during the course of a treatment. Therefore, evaluating on patients temporal data can directly reflect the effectiveness of the $K$-adaptability approach in a clinical context.\par

It is important to note that our developed $K$-adaptability heuristic is a general-purpose method. The heuristic is applicable not only to proton therapy planning but also to other types of robust optimization problems with discrete uncertainty sets.

\section{Conclusion}
In this study, we developed a heuristic for $K$-adaptability problem with discrete uncertainty set. The heuristic outperformed the $K$-medoids method and the LSP and AOSG variants in both worst-case objective estimation and time efficiency. We illustrated the order of solution generation would affect the objective value estimated by the heuristic by comparing it with the AOSG variant. From the clinical perspective, we demonstrated employing the $K$-adaptability approach can improve target coverage significantly compared to the conventional robust optimization approach. The $K$ saturation value evaluated by the $K$-adaptability heuristic indicates that the maximum worst-case objective value can be obtained by a number of treatment plans that is much smaller than the size of the uncertainty set. This number can serve as an indicator of the number of plans to prepare in order to achieve optimal worst-case performance under limited planning resources. Moreover, the developed heuristic is also applicable to other fields other than proton therapy robust treatment planning.

\section*{Acknowledgment}

We would like to thank Bram Gorissen for providing radiotherapy treatment plan optimizer Nymph (\cite{Nymph}) for this study. 

\section*{Funding}
This project has received funding from the European Union’s Horizon 2020 Marie Sklodowska-Curie Actions under Grant Agreement No.955956. It has also been supported in part by grants R01CA266275 and R01CA266803 from the US National Cancer Institute. The content is solely the responsibility of the authors and does not necessarily represent the official views of the National Institutes of Health.




\appendix
\renewcommand{\thesection}{\Alph{section}} 
\renewcommand{\theequation}{\thesection.\arabic{equation}} 
\renewcommand{\thealgorithm}{\thesection.\arabic{algorithm}}
\renewcommand{\thefigure}{\thesection.\arabic{figure}} 

\section{Local solution pool heuristic variant}

The local solution pool heuristic variant (LSP) differs from the proposed $K$-adaptability heuristic in the initial solution pool in the \textit{solution generation} phase. In the LSP heuristic, the initial solution pool starts with the nominal solutions generated for each uncertainty scenario at each $K$ value. In contrast, the initial solution pool of the $K$-adaptability heuristic includes the nominal solutions and the solutions generated for previous $K$ values.

\begin{algorithm}[H]

\setcounter{algorithm}{0}
\caption{Local solution pool (LSP) heuristic variant}\label{alg:LSP_heuristic_variant}
\vskip6pt
\textbf{Input:} Uncertainty scenario set $\mathcal{U}$ = $\{\boldsymbol{D}^1, \dots, \boldsymbol{D}^{T}\}$.\\
\textbf{Output:} Best found combination of $K$ solutions with respect to $\mathcal{U}$ in terms of worst-case objective value $w^*_K$: $X_{K}^*$ = $\{\boldsymbol{x}_1^{(K)}, \dots, \boldsymbol{x}_K^{(K)}\}$, for $K$ = 1,\dots,$|\mathcal{U}|$.
\begin{algorithmic}[1]
\Algphase{Phase 1 - Solution generation}
\State Initialize $K \gets |\mathcal{U}|$ and global solution pool $\mathcal{X}_{global} \gets \emptyset$.
\State Generate an optimal solution for each scenario in $\mathcal{U}$ to form initialization solution set $\mathcal{X}_{int}$.
\While{$ K\ge1$}
    
    \State Initialize scenario clusters pool $\boldsymbol{\mathcal{P}} \gets \emptyset$, $\boldsymbol{C}_{K} \gets \emptyset$, and inner solution pool $\mathcal{X}_{local}\gets\mathcal{X}_{int}$.
    \Repeat
        \State $\boldsymbol{\mathcal{P}}\gets \boldsymbol{\mathcal{P}} \cup \boldsymbol{C}_{K}$
        \State Solve problem (\ref{eq: WC_FL}) on  $\mathcal{X}_{local}$. Obtain   $w^{*}$.
        \State Solve problem (\ref{eq: avg_FL}) on  $\mathcal{X}_{local}$ subject to $w^{*}$. Obtain $\boldsymbol{C}_{K}$. 
         \State Solve problem (\ref{eq:robust planning problem}) for each scenario cluster in $\boldsymbol{C}_{K}$. Obtain solutions   $\mathcal{X}_{R}$.
        \State  $\mathcal{X}_{local} \gets  \mathcal{X}_{local} \cup \mathcal{X}_{R}$.
    \Until $\boldsymbol{C}_{K} \in \boldsymbol{\mathcal{P}}$
    \State$\mathcal{X}_{global} \gets\mathcal{X}_{global} \cup \mathcal{X}_{local}$
    \State $K \gets  K - 1$ .

\EndWhile 
\Algphase{Phase 2 - Solution re-distribution}
\State Set $K = |\mathcal{U}|$.
\While{$ K\ge1$}
    \State Solve problem (\ref{eq: WC_FL}) on $\mathcal{X}_{global}$. Obtain $w^*_K$.
    \State Solve problem (\ref{eq: avg_FL}) on $\mathcal{X}_{global}$ subject to $w^*_K$. Obtain $X_{K}^*$.
    \State $K \gets  K - 1$.

\EndWhile
\end{algorithmic}
\end{algorithm}

\section{Ascending order solution generation heuristic variant}

The ascending order solution generation (AOSG) variant differs from the $K$-adaptability heuristic in the order of $K$ for which the solutions are generated in the \textit{solution generation} phase. The AOSG variant generates solutions for $K$ in an ascending order, namely from $1$ to $|\mathcal{U}|$. In contrast, the $K$-adaptability heuristic generates solution in a descending order in $K$.

\setcounter{algorithm}{0}

\begin{algorithm}[H]
\caption{Ascending order solution generation (AOSG) heuristic variant}\label{alg:AOSG_heuristic_variant}

\vskip6pt
\textbf{Input:} Uncertainty scenario set $\mathcal{U}$ = $\{\boldsymbol{D}^1, \dots, \boldsymbol{D}^{T}\}$.\\
\textbf{Output:} Best found combination of $K$ solutions with respect to $\mathcal{U}$ in terms of worst-case objective value $w^*_K$: $X_{K}^*$ = $\{\boldsymbol{x}_1^{(K)}, \dots, \boldsymbol{x}_K^{(K)}\}$, for $K$ = 1,\dots,$|\mathcal{U}|$.
\begin{algorithmic}[1]
\Algphase{Phase 1 - Solution generation}
\State Initialize $K \gets 1$ and global solution pool $\mathcal{X}_{global} \gets \emptyset$.
\State Generate an optimal solution for each scenario in $\mathcal{U}$ to form initialization solution set $\mathcal{X}_{int}$.
\State $\mathcal{X}_{global} \gets  \mathcal{X}_{global} \cup \mathcal{X}_{int}$ .
\While{$ K\ge1$}
    \State Initialize scenario clusters pool $\boldsymbol{\mathcal{P}} \gets \emptyset$ and $\boldsymbol{C}_{K} \gets \emptyset$
    \Repeat
        \State $\boldsymbol{\mathcal{P}}\gets \boldsymbol{\mathcal{P}} \cup \boldsymbol{C}_{K}$
        \State Solve problem (\ref{eq: WC_FL}) on  $\mathcal{X}_{global}$. Obtain   $w^{*}$.
        \State Solve problem (\ref{eq: avg_FL}) on  $\mathcal{X}_{global}$ subject to $w^{*}$. Obtain $\boldsymbol{C}_{K}$. 
         \State Solve problem (\ref{eq:robust planning problem}) for each scenario cluster in $\boldsymbol{C}_{K}$. Obtain solutions   $\mathcal{X}_{R}$.
        \State  $\mathcal{X}_{global} \gets  \mathcal{X}_{global} \cup \mathcal{X}_{R}$.
    \Until $\boldsymbol{C}_{K} \in \boldsymbol{\mathcal{P}}$
   \State $K \gets  K + 1$ 
\EndWhile 
\Algphase{Phase 2 - Solution re-distribution}
\State Set $K \gets 1$.
\While{$ K\ge1$}
    \State Solve problem (\ref{eq: WC_FL}) on  $\mathcal{X}_{global}$. Obtain $w^*_K$.
    \State Solve problem (\ref{eq: avg_FL}) on  $\mathcal{X}_{global}$ subject to $w^*_K$. Obtain $X_{K}^*$.
    \State $K \gets  K +1$
\EndWhile
\end{algorithmic}
\end{algorithm}

\section{$K$-medoids}\label{Appendix:k-medoids}
\setcounter{equation}{0} 
We quantify the distance between two dose influence matrices by the $L$-2 norm: 
\begin{equation*}
   L_{ab} = ||\boldsymbol{D_a}-\boldsymbol{D_b}||_2.
\end{equation*} We use the following MIP problem to solve the $K$-medoids problem exactly:

\begin{equation}\label{eq: $K$-medoid}
\begin{aligned}
\min_{\boldsymbol{y},\boldsymbol{z}} \quad &\sum_{a=1}^{|\mathcal{U}|}\sum_{b=1}^{|\mathcal{U}|}L_{ab}z_{ab}\\
\textrm{s.t.} \quad 
&\sum_{a=1}^{|\mathcal{U}|}y_a \leq K,  \\
&\sum_{a=1}^{|\mathcal{U}|}z_{ab} = 1, \quad b = 1,\dots,|\mathcal{U}|,\\ 
  &z_{ab}, \leq y_a, \quad a = 1,\dots,|\mathcal{U}|, \quad b = 1,\dots,|\mathcal{U}|,\\ 
  & y_a = \{0,1\}, \quad a = 1,\dots,|\mathcal{U}|,\\
  &z_{ab} = \{0,1\}, \quad a = 1,\dots,|\mathcal{U}|, \quad b = 1,\dots,|\mathcal{U}|.\\
\end{aligned}
\end{equation}\par

\section{Results on individual cases}\label{Appendix:individual result}
This section reports the worst-case dose distribution and objective values comparison of case one, two, three, and five.

\setcounter{figure}{0} 

\begin{figure}[H]
    
    \centering
    \caption{Worst-case dose distribution of case one.}
    \includegraphics[width=1.0\textwidth]{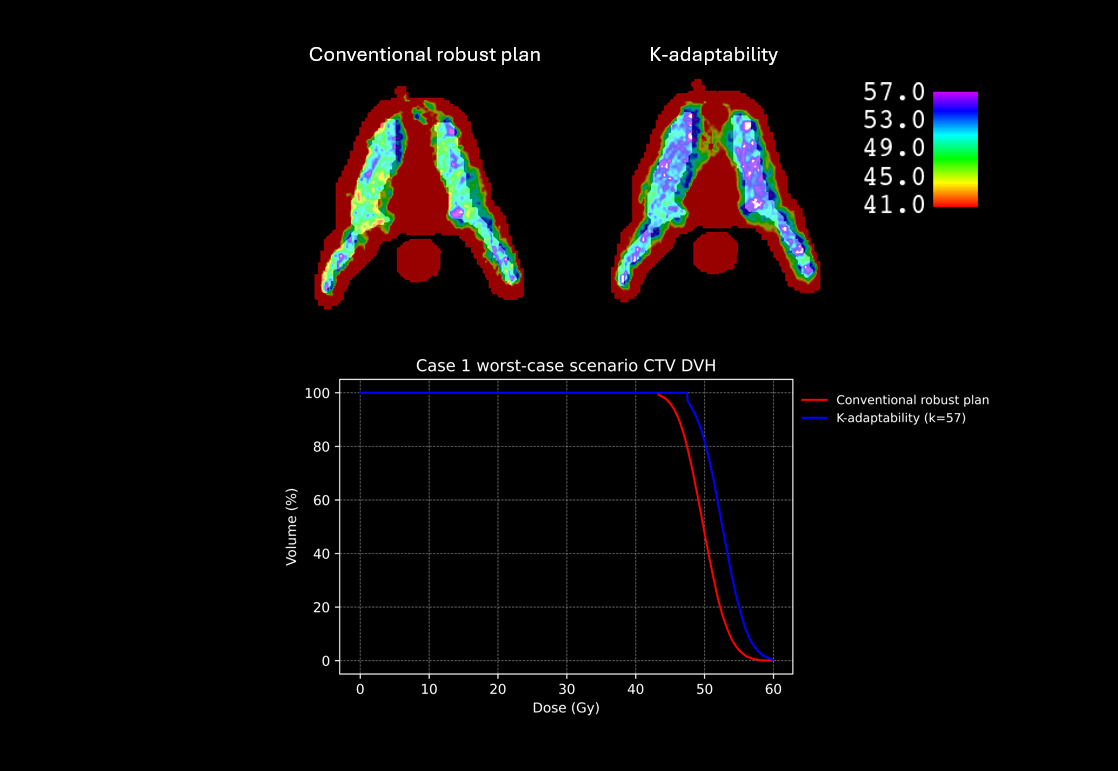}
    \label{fig:case_1_dvh}

\end{figure}

\begin{figure}[H]
    
    \centering
    \caption{Worst-case CTV $D_{min}$ vs $K$ plot of case one.}
    \includegraphics[width=1.0\textwidth]{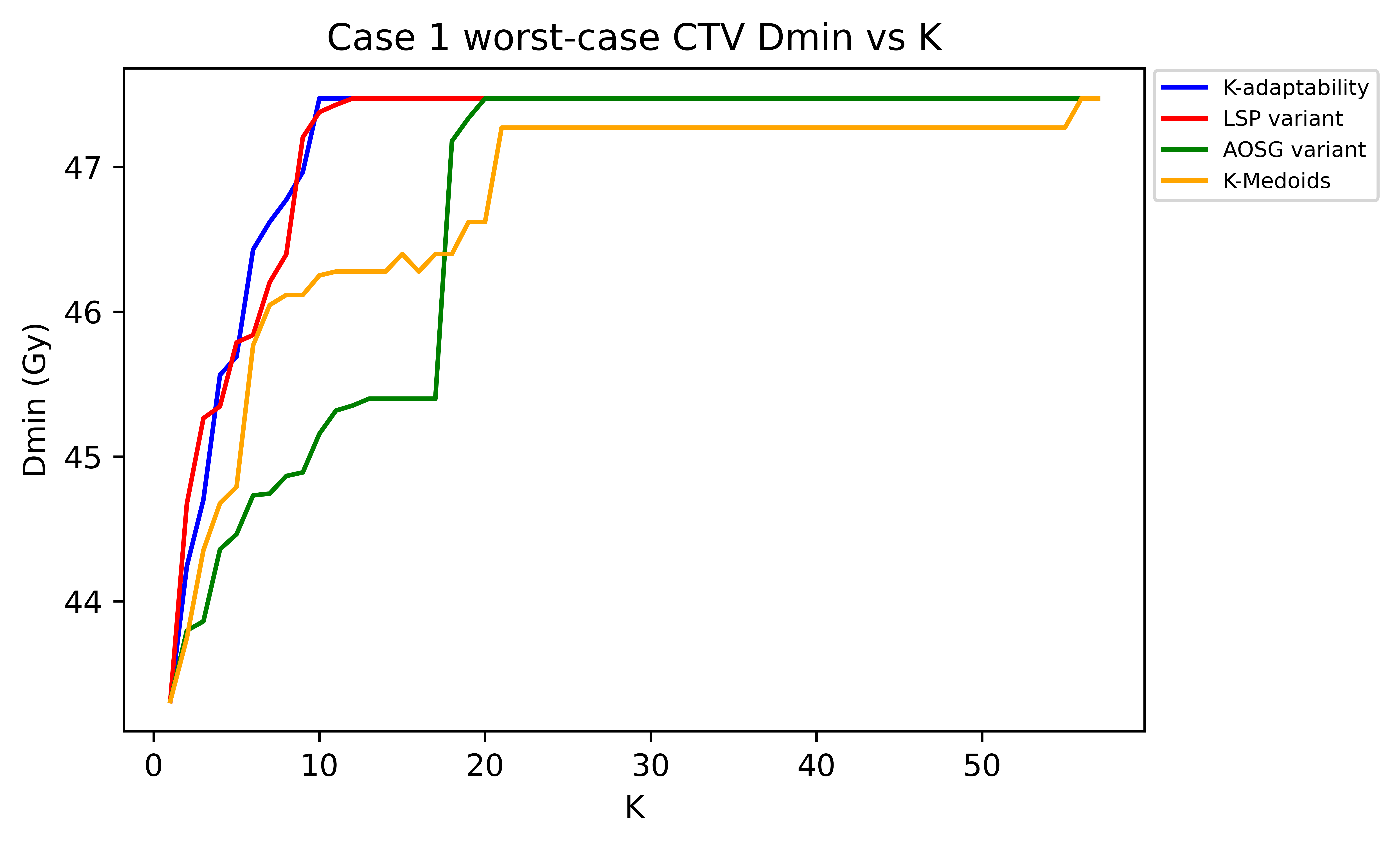}
    \label{fig:case_1_worst_case_objective_vs_K}

\end{figure}

\begin{figure}[H]
    \centering
    \caption{Worst-case dose distribution of case two.}
    \includegraphics[width=1.0\textwidth]{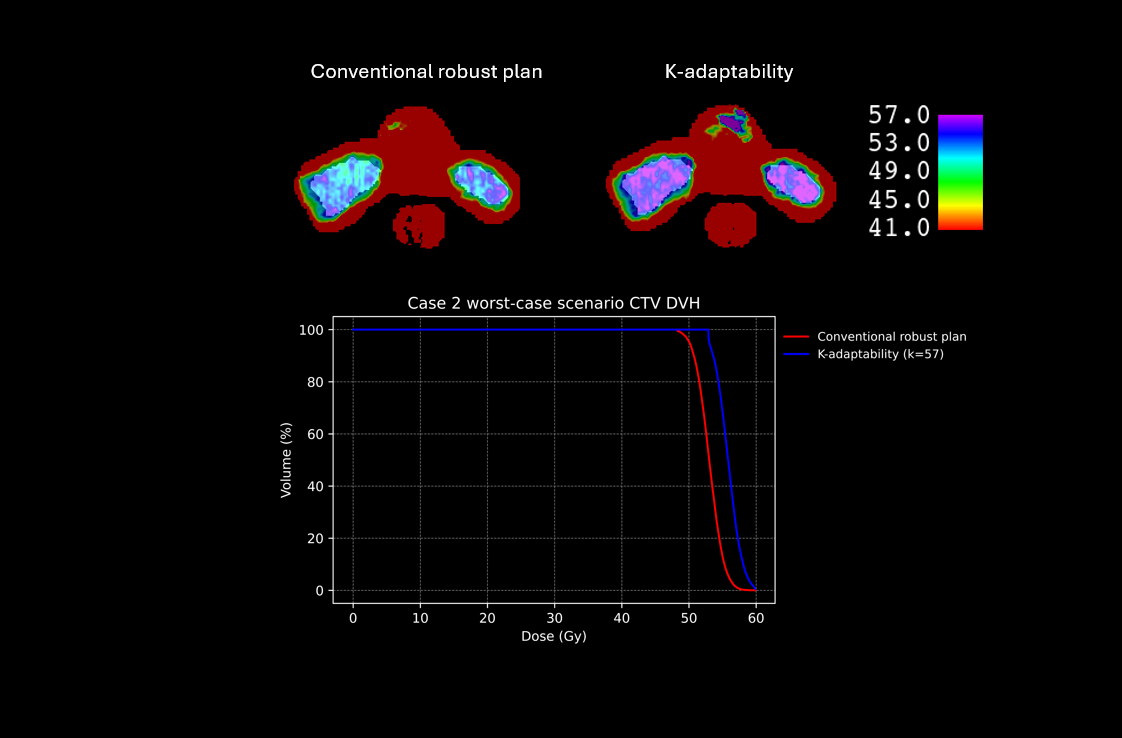}
    \label{fig:case_2_dvh}
\end{figure}

\begin{figure}[H]
    \centering
    \caption{Worst-case CTV $D_{min}$ vs $K$ plot of case two.}
    \includegraphics[width=1.0\textwidth]{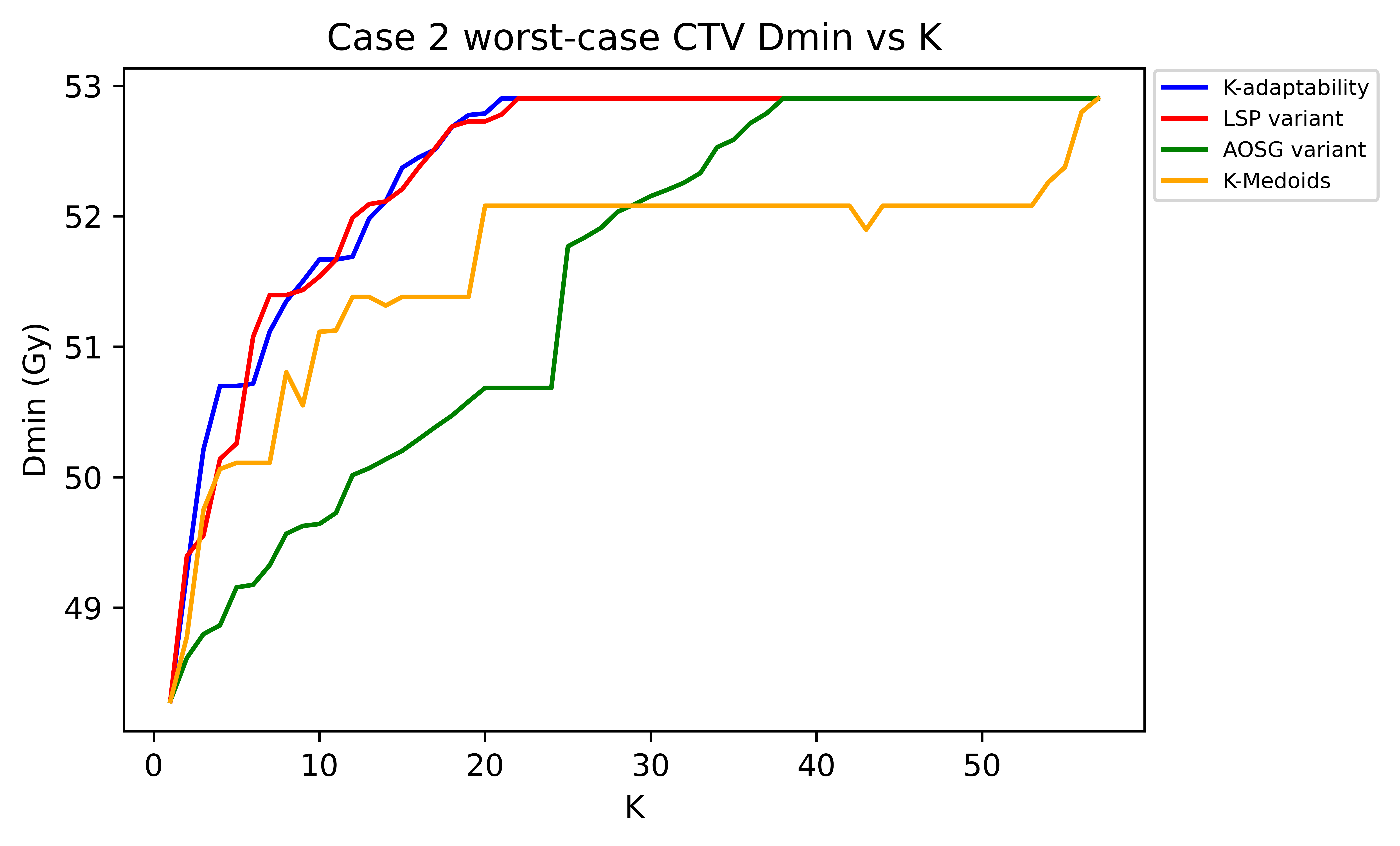}
    \label{fig:case_2_worst_case_objective_vs_K}
\end{figure}

\begin{figure}[H]
    \centering
    \caption{Worst-case dose distribution of case three.}
    \includegraphics[width=1.0\textwidth]{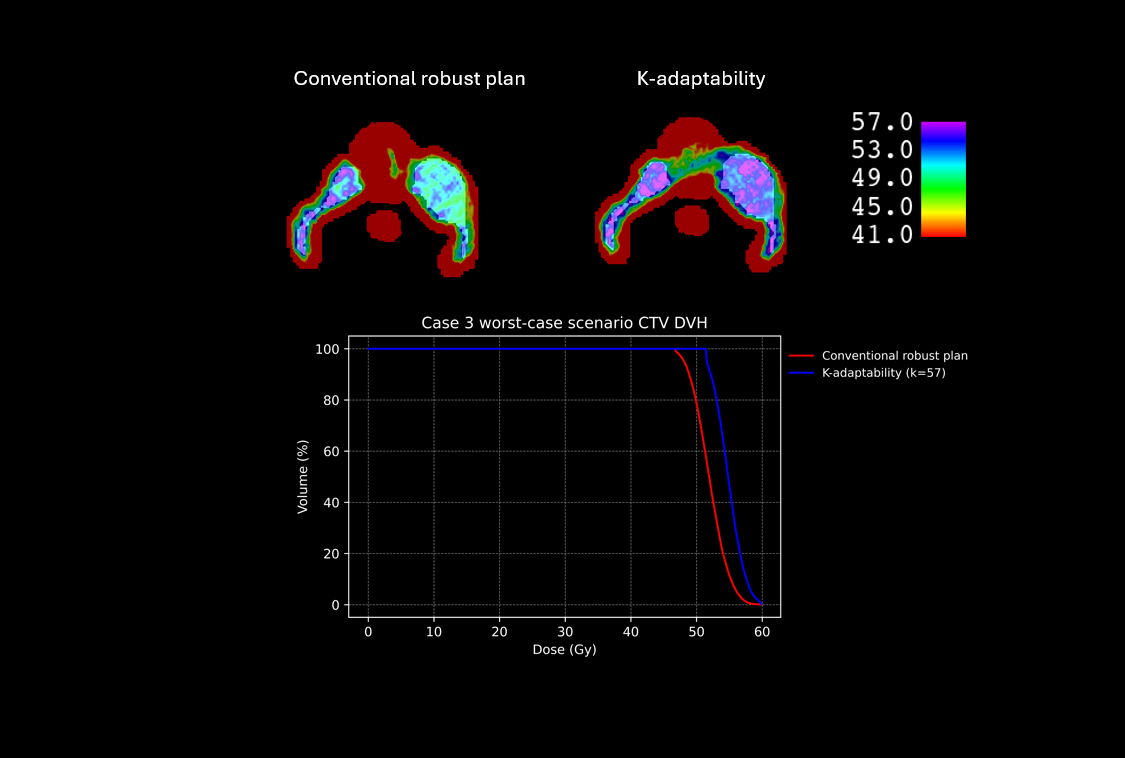}
    \label{fig:case_3_dvh}
\end{figure}

\begin{figure}[H]
    \centering
    \caption{Worst-case CTV $D_{min}$ vs $K$ plot of case three.}
    \includegraphics[width=1.0\textwidth]{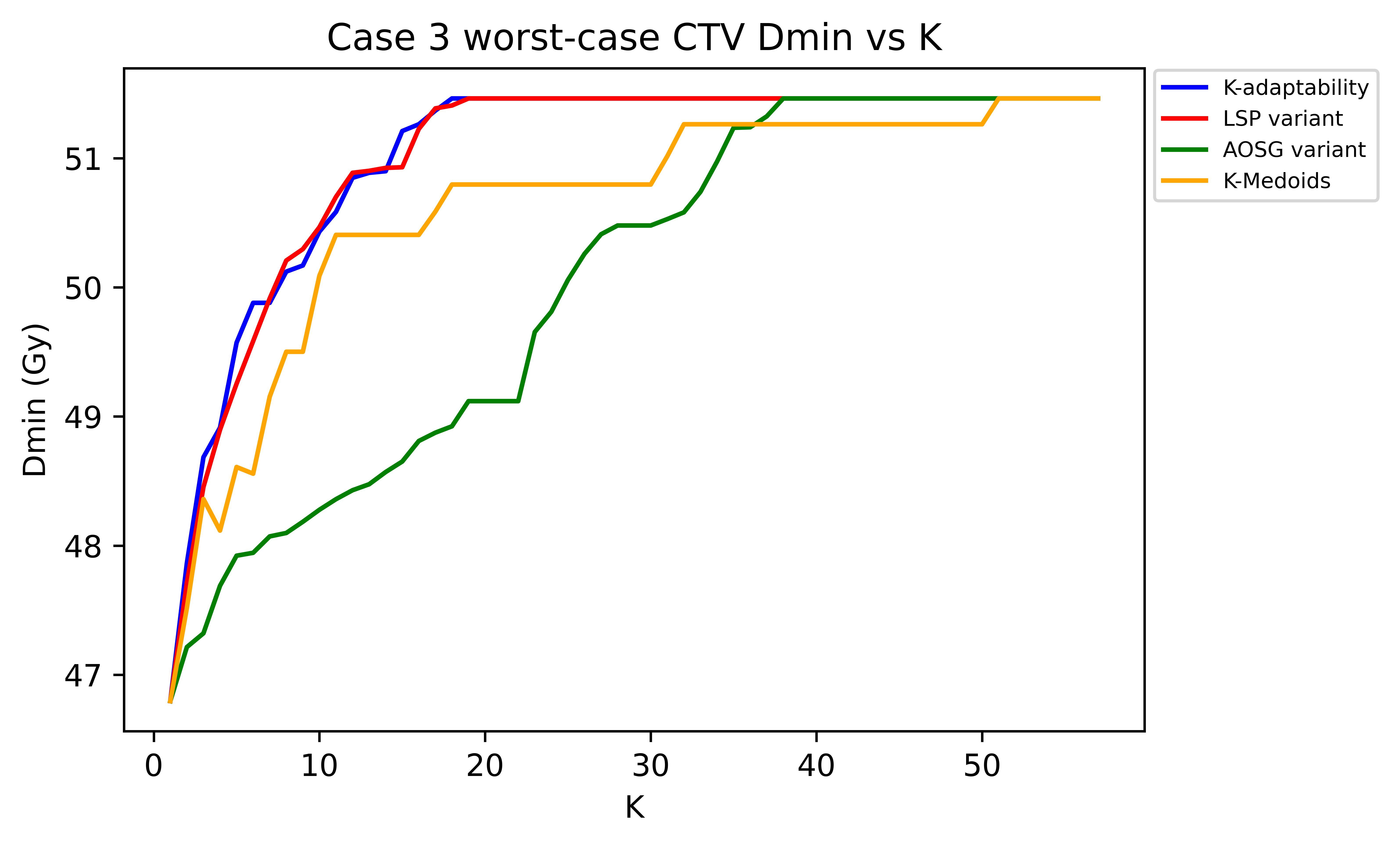}
    \label{fig:case_3_worst_case_objective_vs_K}
\end{figure}

\begin{figure}[H]
    \centering
    \caption{Worst-case dose distribution of case five.}
    \includegraphics[width=1.0\textwidth]{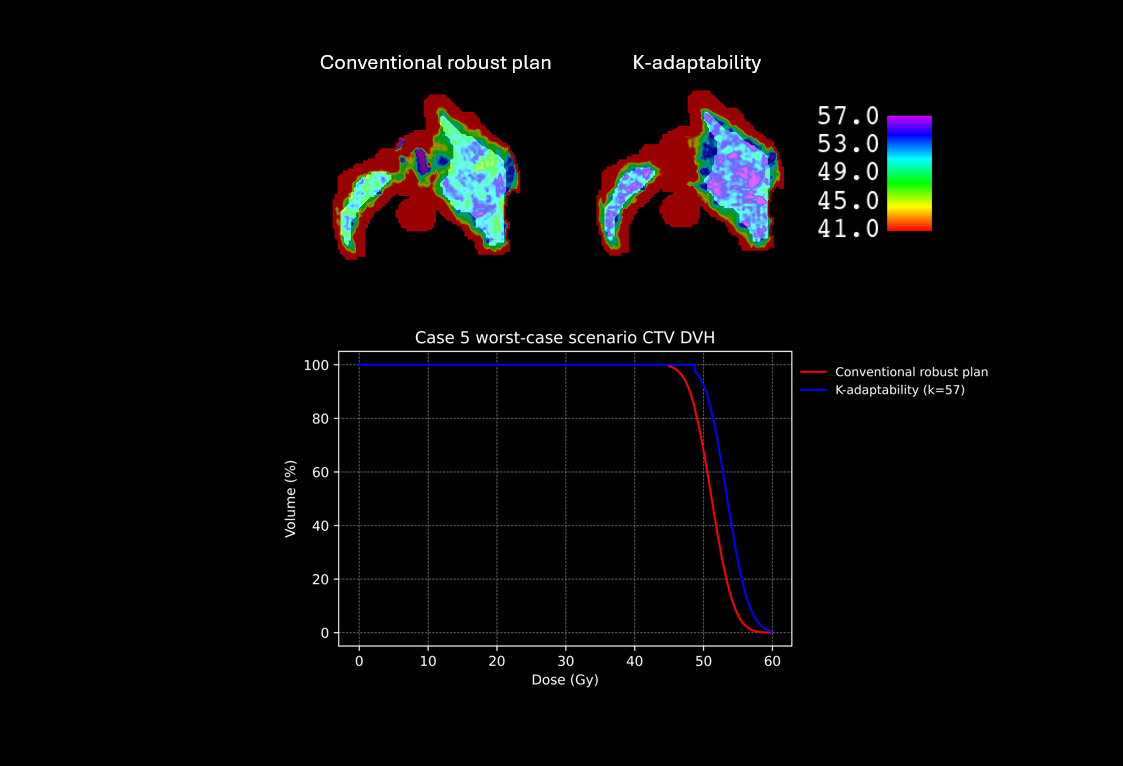}
    \label{fig:case_5_dvh}
\end{figure}

\begin{figure}[H]
    \centering
    \caption{Worst-case CTV $D_{min}$ vs $K$ plot of case five.}
    \includegraphics[width=1.0\textwidth]{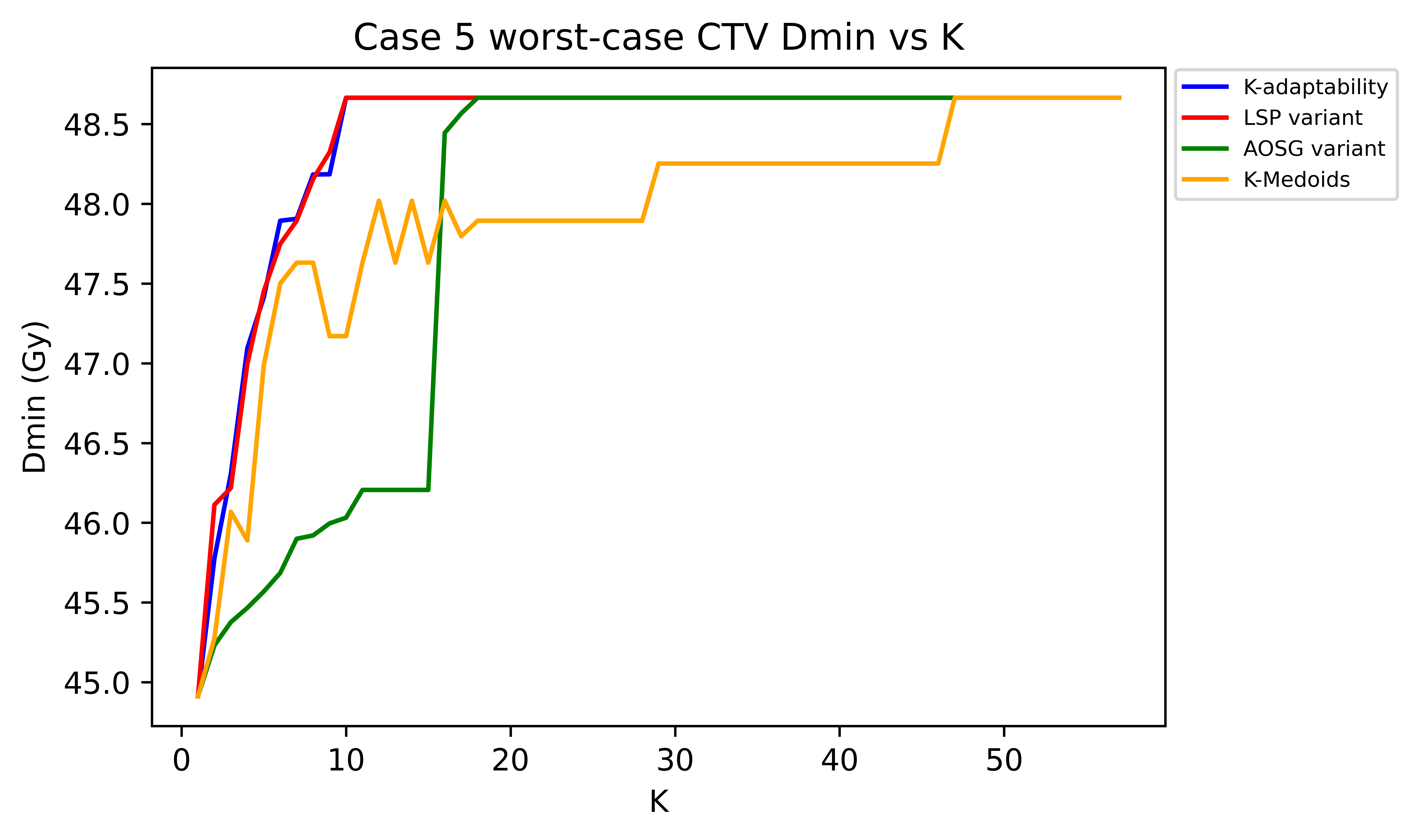}
    \label{fig:case_5_worst_case_objective_vs_K}
\end{figure}

\bibliographystyle{apalike}

\end{document}